\documentclass[preprint,12pt]{elsarticle}

\usepackage{mathptmx}       
\usepackage{helvet}         
\usepackage{courier}        
\usepackage{type1cm}        
\usepackage{makeidx}         
\usepackage{graphicx}        
\usepackage{multicol}        
\usepackage[bottom]{footmisc}
\usepackage{amsmath,amssymb,amscd}
\usepackage{amsfonts}
\usepackage{enumitem}

\usepackage[
colorlinks=true,
citecolor=blue,
urlcolor=blue,
linktocpage,
pagebackref
]{hyperref}

\usepackage{comment}

\def\bx{{\bf x}}   
\def\cmb{{\omega}}   
\def\lim{{\xi}}    

\journal{Computers \& Fluids}

\begin{document}

\begin{frontmatter}



\title{Simulation-Driven Optimization of High-Order Meshes
       in ALE Hydrodynamics \tnoteref{l_title}}

\author{Veselin Dobrev \fnref{llnl}}
\author{Patrick Knupp \fnref{dihedral}}
\author{Tzanio Kolev \fnref{llnl}}
\author{Ketan Mittal \fnref{llnl}}
\author{Robert Rieben \fnref{llnl}}
\author{Vladimir Tomov \corref{cor1} \fnref{llnl}}
\fntext[llnl]
{Lawrence Livermore National Laboratory, 7000 East Avenue, Livermore, CA 94550}
\fntext[dihedral]
{Dihedral LLC, Bozeman, MT 59715}
\cortext[cor1]
{Corresponding author, tomov2@llnl.gov}
\tnotetext[l_title]
{Performed under the auspices of the U.S. Department of Energy under
Contract DE-AC52-07NA27344 (LLNL-JRNL-800087)}

\address{}

\begin{abstract}

In this paper we propose tools for high-order mesh optimization and demonstrate their benefits in the
context of multi-material Arbitrary Lagrangian-Eulerian (ALE) compressible shock hydrodynamic applications.
The mesh optimization process is driven by information provided by the
simulation which uses the optimized mesh, such as shock positions,
material regions, known error estimates, etc.
These simulation features are usually represented discretely,
for instance, as finite element functions on the Lagrangian mesh.
The discrete nature of the input is critical for the practical applicability of the algorithms we
propose and distinguishes this work from approaches that strictly require analytical
information.
Our methods are based on node movement through a high-order extension of the
Target-Matrix Optimization Paradigm (TMOP) of \cite{Knupp2012}.
The proposed formulation is fully algebraic and relies
only on local Jacobian matrices, so it is applicable to all types of mesh elements,
in 2D and 3D, and any order of the mesh.
We discuss the notions of constructing adaptive target matrices and obtaining
their derivatives, reconstructing discrete data in intermediate meshes,
node limiting that enables improvement of global mesh quality while
preserving space-dependent local mesh features, and appropriate normalization of the
objective function.
The adaptivity methods are combined with automatic ALE triggers that can
provide robustness of the mesh evolution and avoid excessive remap procedures.
The benefits of the new high-order TMOP technology are illustrated on several
simulations performed in the high-order ALE application BLAST \cite{blast}.
\end{abstract}

\end{frontmatter}


\section{Introduction}
\label{sec_intro}

Lagrangian methods for compressible multi-material shock hydrodynamics
\cite{VonNeumann1950, Benson1992, Loubere2004, Scovazzi2007, Maire2009} are
characterized by a computational mesh that moves with the material velocity.
A key advantage of using such methods is that they do not produce any
numerical dissipation around contact discontinuities and material interfaces,
when these are aligned with the element boundaries of the moving mesh.
This is a consequence of the fact that the dynamics equations, when written in
a Lagrangian frame, do not include nonlinear convective terms.
However, the main disadvantage of Lagrangian methods is that they lead to
mesh distortion, which results in small time steps or simulation breakdowns.
Many problems of interest, i.e., impact simulations, cannot be run to
longer evolution times by a purely Lagrangian method.
These disadvantages are addressed by the Arbitrary Lagrangian-Eulerian (ALE)
approach \cite{Hirt1974, Barlow2016, Maire2007, Galera2010}, which involves the notion of
periodic mesh optimization.
However, while improving the mesh, ALE methods can produce
their own numerical errors around the material interfaces, thus counteracting
the main advantage of the Lagrangian approach.

In this paper we propose mesh optimization methods that provide mechanisms
to control the amount of numerical error resulting from the ALE approach.
Keeping unchanged the method's underlying solution transfer (remap) and
Lagrangian phases, we explore how mesh adaptivity through node movement
(a.k.a. $r$-adaptivity), adaptive ALE triggers,
and restricted node movement can be used to decrease
the errors resulting from the overall ALE procedure.
This paper extends our previous work on ALE hydrodynamics in the BLAST code
\cite{Dobrev2012,Dobrev2016,Dobrev2018}
and mesh optimization by the Target-Matrix Optimization Paradigm
\cite{Knupp2012, ETHOS2018}.
We incorporate the TMOP technology in the context of a high-order finite element
ALE approach for multimaterial flow and introduce several important
notions, including limiting of node movement based on the simulation dynamics,
proper normalization of the energy functional, adaptivity to discrete
simulation features, and adaptive ALE triggers.
We focus on high-order curved meshes which have recently shown several
mathematical and computational advantages, including optimal convergence rates
on domains with curved boundaries/interfaces, symmetry preservation in radial
flow \cite{Dobrev2012, Dobrev2013} and
equivalent simulation quality with a smaller number of degrees of freedom
\cite{Dobrev2012, Tomov2016, Boscheri2016}.

Although in this paper we consider a particular application (BLAST), an
important aspect of the TMOP approach is its generality.
In particular, our methods are fully algebraic and do not involve any
low-level geometric manipulations.
They only require information about the local Jacobians, and thus can be applied to
meshes of any order (including linear meshes), element type, and dimension.
The proposed methods can be applied to both monolithic ALE methods
\cite{Scovazzi2016, Bakosi2017, Gaburro2017} to define their mesh velocity,
and ALE methods that split the Lagrangian, mesh optimization and remap phases.
They can also be applied to diffused interface methods, where ALE errors are
decreased by shrinking the transition region of the volume fraction functions,
or exact interface representation methods
\cite{Kucharik2010}, where ALE errors can be decreased by obtaining better
mesh resolution in the interface regions.

The algorithms presented in this paper rely strictly on node movement and
preserve the original topology of the mesh.
Reconnection-based remesh approaches are more powerful in terms of
controlling the mesh characteristics, but also require more complicated
solution transfer procedures, see \cite{Loubere2010}.

Our mesh optimization method is based on the TMOP framework described in
\cite{ETHOS2018, IMR2018, knupp2019target}.
This approach is distinguished from similar methods by its emphasis on
target-matrix construction methods that permit a greater degree of control
over the optimized mesh.
Pointwise mesh quality metrics are defined by utilizing sub-zonal information
obtained by sampling element Jacobians at element quadrature points.
These metrics are capable of measuring shape, size or alignment of
the region around the point of interest.
TMOP requires predefined target-matrices as a way for the user to incorporate
application-specific physical information into the metric that is being
optimized.
The combination of targets and quality metrics is used to optimize the node
positions, so that they are as close as possible to the shape, size and/or alignment of
their targets.

The rest of the paper is organized as follows.
In Section \ref{sec_prelim} we review the multi-material ALE hydrodynamics framework
in BLAST and the basic TMOP components.
Section \ref{sec_opt} describes the TMOP optimization functional in
further detail, along with the concepts of limited node motion and
proper normalization of the resulting nonlinear functional.
Mesh adaptivity to discrete simulation features and its relation to
dynamic ALE triggers is discussed in Section \ref{sec_adapt}.
Section \ref{sec_tests} presents several numerical tests that demonstrate the
features of the discussed methods.
Conclusions are presented in Section \ref{sec_concl}.



\section{Preliminaries}
\label{sec_prelim}

In this section we provide a brief overview of the multi-material ALE framework
(as implemented in the BLAST code) that is used to test the proposed mesh optimization methods.
We focus only on aspects that are relevant to our mesh optimization
purposes, while a complete description of our ALE hydrodynamics
methods can be found in \cite{Dobrev2018}.
We also provide a summary of the basic TMOP components in the context of BLAST.


\subsection{Multi-Material ALE Framework}
\label{sec_blast}

BLAST solves the multi-material Euler equations by clearly separating the
Lagrange, remesh and remap phases.
It utilizes a single-fluid multi-material description of the system where
all materials share a common velocity, but each has its own density and energy.
BLAST supports performing an arbitrary number of Lagrangian steps between two
ALE steps (remesh+remap), as opposed to \emph{continuous remap} algorithms
that perform an ALE step after every Lagrangian step.

\paragraph{Space discretization}
The space discretization is based on high-order finite elements.
Velocity $v$ is discretized in the finite element space
$\mathcal{V} \subset [H^1(\widehat{\Omega})]^d$, with basis $\{ w_i \}$,
where $\widehat{\Omega}$ is the computational mesh corresponding to the
domain of interest $\Omega$, and $d$ is the space dimension.
Specific internal energies $\{e_m\}$ are discretized in
$\mathcal{E} \subset L^2(\widehat{\Omega})$, with basis $\{ \phi_j \}$, where
$m$ is the material index.
Material densities $\{\rho_m\}$ are evolved at certain quadrature points of
interest through the notion of pointwise mass conservation.
Throughout this manuscript we refer to the pairs of spaces $Q_k Q_{k-1}$,
by which we denote $\mathcal{V} = (Q_k)^d$, i.e., the Cartesian product of
the space of continuous finite elements on quadrilateral or hexahedral
meshes of degree $k$, and $\mathcal{E} = \hat{Q}_{k-1}$, the companion space of
discontinuous finite elements of order one less than the kinematic space.
As shown in \cite{Dobrev2012} the specific choice of $Q_1 Q_0$ corresponds
to traditional staggered grid discretizations under additional simplifying
assumptions.

\paragraph{Discrete mesh representation}
To obtain a discrete representation of the high-order mesh,
we start with the set of kinematic scalar basis functions
$\{ \bar{w}_i \}_{i=1}^{N_w}$ on a reference element $\bar{E}$.
The shape of any element $E$ in the mesh is then fully described by a
matrix $\mathbf{x}_E$ of size $d \times N_w$ whose columns represent the
coordinates of the element control points (a.k.a. nodes or element degrees of freedom).
Given $\mathbf{x}_E$, we introduce the map $\Phi_E:\bar{E}
\to \mathbb{R}^d$ whose image is the high-order element $E$:
\begin{equation}
\label{eq_x}
x(\bar{x}) =
   \Phi_E(\bar{x}) \equiv
   \sum_{i=1}^{N_w} \mathbf{x}_{E,i} \bar{w}_i(\bar{x})\,,
   \qquad \bar{x} \in \bar{E},
\end{equation}
where we used $\mathbf{x}_{E,i}$ to denote the $i$-th column of $\mathbf{x}_E$.
To ensure continuity between mesh elements, we define a global vector
$\mathbf{x}$ of mesh positions that contains the $\mathbf{x}_E$ control points
for every element.
For any element $E$ in the mesh, we can compute the Jacobian of the mapping
$\Phi_E$ at any reference point $\bar{x} \in \bar{E}$ as
\begin{equation}
\label{eq_A}
  A_E(\bar{x}) = \frac{\partial \Phi_E}{\partial \bar{x}} =
    \sum_{i=1}^{N_w} \mathbf{x}_{E,i} [ \nabla \bar{w}_i(\bar{x}) ]^T \,.
\end{equation}
We will refer to the initial mesh, which is to be optimized,
as the {\em Lagrangian} mesh.
We denote the Lagrangian mesh by $\mathcal{M}_0$.

\paragraph{Material representation and evolution}
Material positions are represented on a discrete
level as point-wise volume fractions at certain quadrature points of interest.
At any such quadrature point $x$, the volume fraction values $\{\eta_m\}$ satisfy $\eta_m(x) \geq 0$ and $\sum_m \eta_m(x) = 1$.
During the Lagrange phase, these volume fraction are evolved by a
sub-zonal closure model, see \cite{Dobrev2016}.
Before each remap procedure, these quadrature point values are projected to
discontinuous $L^2$ functions, which are remapped,
by solving an advection-based PDE, to the new mesh.
Once the remap step is complete, the remapped $L^2$ functions are interpolated
back to the quadrature points of interest.
The above ALE treatment of the material positions leads to numerical diffusion
around the material interfaces, which is often the main source of numerical
dissipation error in BLAST.


\subsection{Overview of TMOP}
\label{sec_tmopoverview}

Mesh optimization in TMOP is driven by moving the mesh nodes
(i.e. changing the values of the vector $\mathbf{x}$) in a way that minimizes
a mesh quality metric, $\mu(T)$, over the computational mesh.
In this section we explain the meaning of $\mu$ and $T$.
A more detailed discussion about the TMOP optimization functional is
presented in Section \ref{sec_opt}.

The \emph{weighted Jacobian} matrix $T$ represents the transformation from
a \emph{target} configuration to the current \emph{physical} (or \emph{active})
mesh positions, so that a perfect mesh would result in $T = I$.
At each quadrature point, the matrix $T$ is determined using two
Jacobian matrices:
\begin{itemize}
  \item The \emph{target-matrix} $W_{d \times d}$, which is the Jacobian of the
    transformation from the \emph{reference} to the \emph{target} coordinates.
    This matrix is defined according to a user-specified method prior to
    optimization; it defines the desired local properties of the optimal mesh.
  \item The Jacobian matrix $A_{d \times d}$ of the transformation
    between \emph{reference} and \emph{physical} space coordinates,
    always computed by \eqref{eq_A}.
\end{itemize}
Then we have $T=AW^{-1}$.
Note that both $A$ and $W$ are space-dependent functions of the nodal
coordinates.
As any Jacobian matrix, $W$ can be written as a decomposition of matrices that
represent its main geometric properties:
\begin{equation}
\label{eq_decomp}
W = [\text{volume}][\text{orientation}][\text{skew}][\text{aspect ratio}].
\end{equation}
The concept of target construction allows the user specify some of the above
components, and thus control the desired geometric property
\cite{knupp2019target}.
We utilize this concept further in Section \ref{sec_adapt} to achieve
simulation-based mesh adaptivity.

\begin{figure}[h!]
\centerline{
  \includegraphics[width=0.5\textwidth]{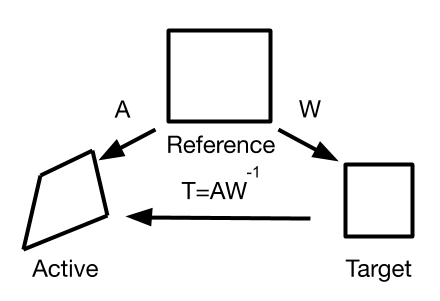}}
\caption{Schematic representation of the major TMOP matrices.}
\label{fig_tmop}
\end{figure}

Given a matrix $T=AW^{-1}$, the purpose of the local quality metric $\mu(T)$ is
to measure the differences between $A$ and $W$ only in one or more particular
components of \eqref{eq_decomp}, while being invariant of the others.
The choice of quality metric is based on the
geometric aspects that the user wants to control.
For example, $\mu_2=\mid T \mid^2/2\tau-1$ is a \emph{shape} metric that
controls skew and aspect ratio, but is invariant to orientation and volume.
Here, $\mid T \mid$ is the Frobenius norm of $T$ and $\tau=\text{det}(T)$.
Similarly, $\mu_7 = \mid T-T^{-t}\mid^2$ and
$\mu_9=\tau\mid T-T^{-t}\mid^2$ are \emph{shape}$+$\emph{size}
metrics that controls volume, skew and aspect ratio,
but are invariant to orientation.
The reader is referred to \cite{ETHOS2018} for a more detailed
discussion on quality metrics.


\section{Recent Improvements in TMOP}
\label{sec_opt}

In this section we go deeper into the details of the TMOP's nonlinear objective
function and present several improvements over our previous formulation
presented in \cite{ETHOS2018, IMR2018}.
These modifications are motivated by practical concerns related to
restricting node movement and proper normalization of the different terms
in the objective function.


\subsection{Objective Function}
\label{sec_objective}

TMOP optimizes the mesh quality by minimizing a global objective
function $F(x)$ that depends on the local quality measure throughout the mesh.
This minimization is performed by solving
$\partial F / \partial \mathbf{x} = 0$ with Newton's method.
The objective function has the following general form:
\begin{equation}
\label{eq_F_full}
  F(x) =
  \frac{1}{n}
    \sum_{s=1}^n
    \frac{\sum_{E(x)} \int_{E_t} \cmb_s(x) \mu_{i_s}(T(x)) dx_t}
         {\sum_{E(x_0)} \int_{E_t} \cmb_s(x_0) \mu_{i_s}(T_0(x_0)) dx_t} +
  c \sum_E \int_{E_t} \xi(x-x_0, \delta(x_0)) dx.
\end{equation}
The right-most term is used to limit the node displacements during optimization.
This term is discussed in Section \ref{sec_limiting}.
The other term represents an explicit combination of $n$ mesh quality metrics
$\mu_{i_1}, \dots \mu_{i_n}$.
Although most problems require only one metric ($n=1$),
in other cases \cite{mittal2019mesh} it might be necessary to combine metrics
with different weights $\cmb$ in order to make some of them more dominant in
certain regions of the domain.
The integrals in \eqref{eq_F_full} are computed as
\begin{equation}
\label{eq_vm}
  \sum_{E \in \mathcal{M}} \int_{E_t} \cmb(x_t) \mu(T(x_t)) dx_t =
  \sum_{E \in \mathcal{M}} \sum_{x_q \in E_t}
                           w_q\,\det(W(\bar{x}_q))\, \cmb(x_q) \mu(T(x_q)),
\end{equation}
where $\mathcal{M}$ is the current mesh,
$E_t$ is the target element corresponding to the physical element $E$,
$w_q$ are the quadrature weights, and the point $x_q$ is the image of the
reference quadrature point $\bar{x}_q$ in the target element $E_t$.
Note that the right-hand side in \eqref{eq_vm} depends on the mesh positions
$x$ through the Jacobian matrices $A$ used in the definition of $T$.
The integration in \eqref{eq_vm} is performed over the target elements,
enforcing that the integral contribution from a given element $E$, relative to
the contributions from other elements, is only based on the difference
with its target $E_t$ (which is measured by $\mu(T)$), and not on its relative
size compared to the other elements.
In particular, very small elements will
not be neglected by the optimization process due to their size.
The existence of a minimum for variational optimization problems
like \eqref{eq_vm} has been explored theoretically in
\cite{Garanzha2014, Garanzha2010}.

In general, the various quality metrics $\mu_i$ have different magnitudes;
we have observed orders of magnitude differences in the resulting integrals.
The role of the denominators in the first sum of \eqref{eq_F_full} is to
normalize the magnitude of the different metrics.
The sum of the metric terms is always equal to $1$ at the start of the
optimization process and is expected to decrease during the iteration.
Furthermore, all metric integrals start at the same value of $1/n$.
This normalization also makes the magnitude of each metric integral
invariant of $i$, the level of refinement,
the specifics of the problem (e.g. the units system),
and the size of the domain.


\subsection{Limiting of Node Displacements}
\label{sec_limiting}

\begin{figure}[h!]
\begin{center}
$\begin{array}{ccc}
\includegraphics[width=65mm]{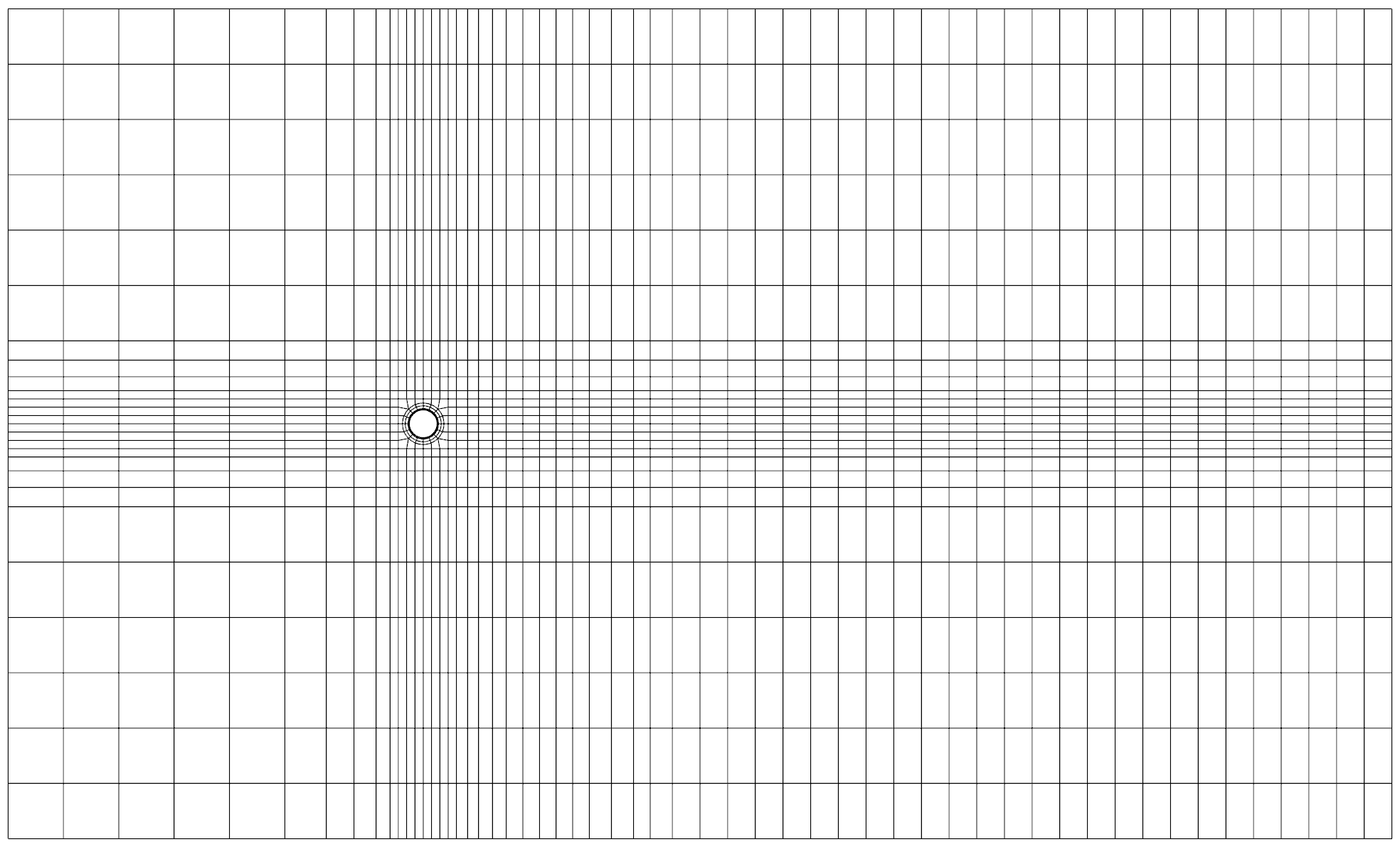} &
\includegraphics[width=65mm]{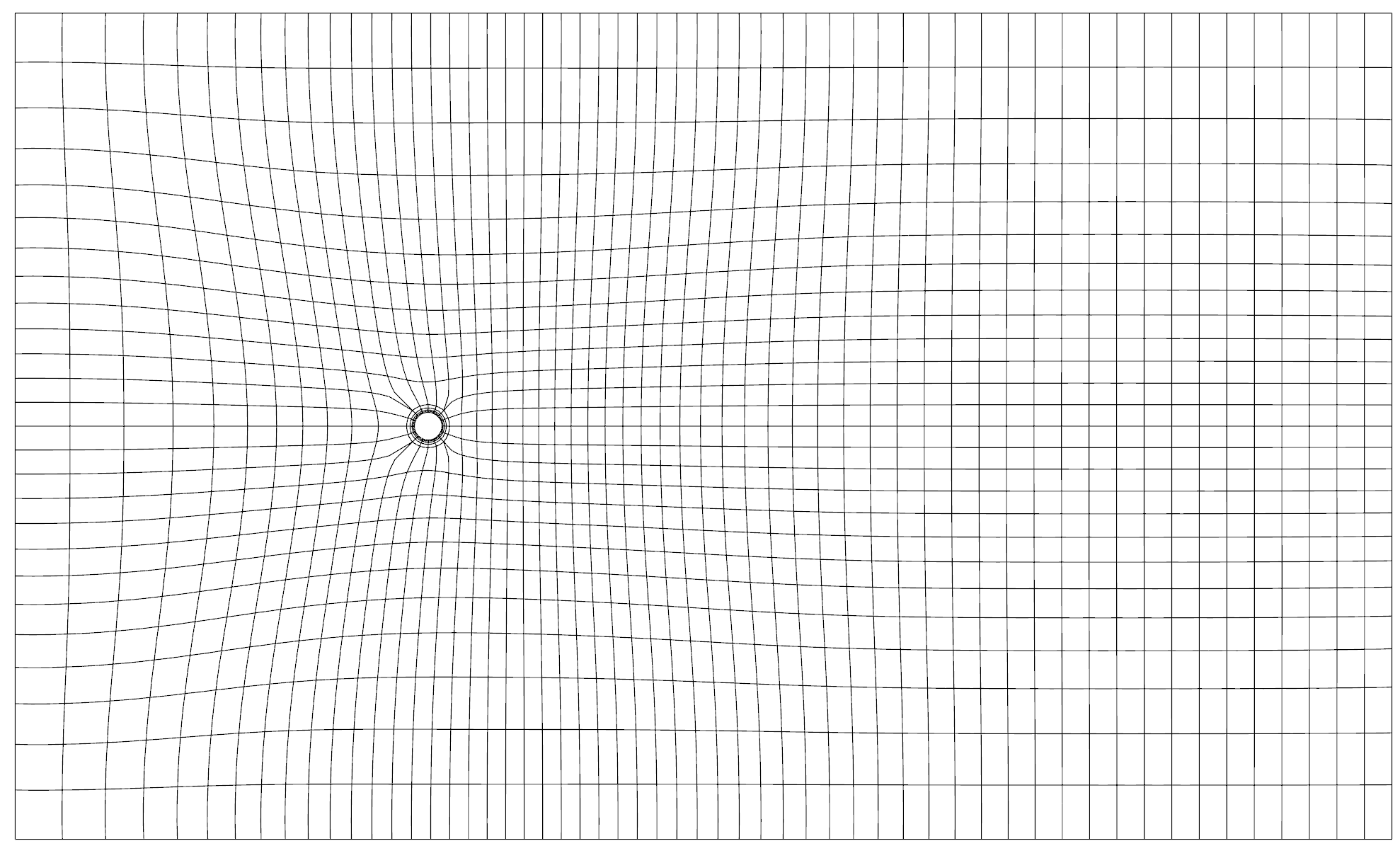} \\
\end{array}$
\end{center}
\vspace{-7mm}
\caption{Original mesh (left) and optimized mesh (right) through
         space-dependent restriction of the node displacements.}
\label{fig_cylinder}
\end{figure}

One advantage of the Lagrangian approach is that
it provides a natural form of mesh adaptivity.
For example, the mesh spacing usually gets compressed around shock regions.
Thus, it might be beneficial to have the ability to preserve, up to some extent,
the positions $x_0$ of the Lagrangian mesh $\mathcal{M}_0$.
Another important benefit of staying close to the Lagrangian mesh is minimizing
errors due to remap.
The purpose of the limiting term (the right-most term in in \eqref{eq_F_full})
is to restrict the movement of nodes during the mesh optimization procedure
with respect to their initial positions $x_0$.
Note that the limiting term is space-dependent, allowing the preservation of mesh
features (e.g. initial resolution) in specific regions of the domain, which
is an important requirement for various applications \cite{mittal2019mesh}.
A simple motivational 2D example is shown in Figure \ref{fig_cylinder}, where
the goal is to preserve the resolution around the central boundary layer,
while achieving uniform mesh resolution in the rest of the domain.

In what follows, we rely on the fact that there exists a
bijective map between every position $x$, in the current mesh $\mathcal{M}$,
and its initial position $x_0$ in the Lagrangian mesh $\mathcal{M}_0$,
which is the case in methods based on node movement.
We define the limiting term from \eqref{eq_F_full} as a general function of
distance between the current positions $x$ and their corresponding Lagrangian
positions $x_0$ as
\begin{equation}
\label{eq_lim}
  c \sum_E \int_{E_t} \xi \left( x-x_0[x], \delta(x_0[x]) \right),
\end{equation}
where $\delta(x_0) > 0$ is a space-dependent scalar function that specifies
the acceptable amount of displacement for each mesh node.
Note that moving a given position $x$ changes the displacement $x-x_0$,
but it does not change its image on the initial mesh $x_0[x]$
and the value of $\delta(x_0[x])$.
The values of $\delta$ are problem-specific and must be specified by the user.
For example, in the remesh phase of ALE simulations, $\delta(x_0)$ can be some
fraction of the displacement that occurs between the last two ALE steps.

\begin{figure}[h!]
\centerline{\includegraphics[width=0.5\textwidth]{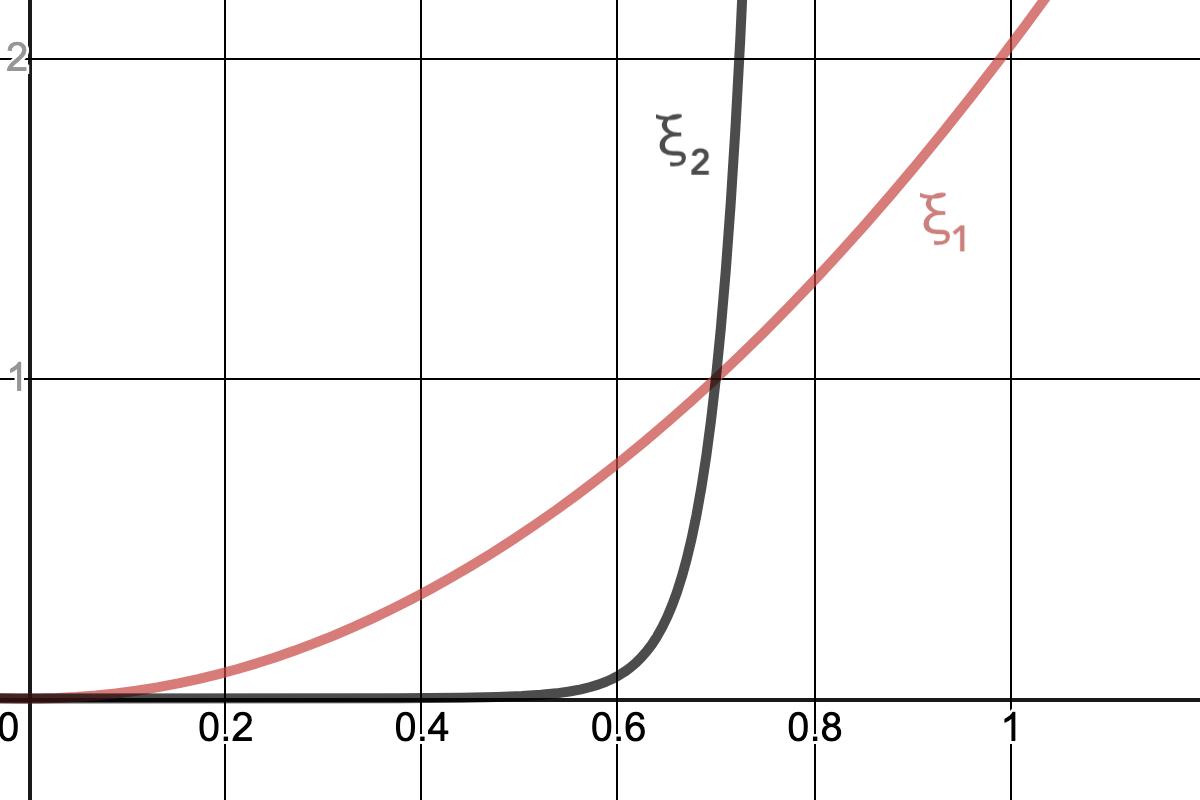}}
\caption{Graphs of $\xi_1$ (red) and $\xi_2$ (black) for $\delta = 0.7$.}
\label{fig_xi}
\end{figure}

To achieve the limiting effect, the function $\xi$ has the following properties:
\begin{enumerate}
\item $\xi \rightarrow 0$ when $x \rightarrow x_0$, and in particular,
      the initial value of the limiting integral is negligible.
\item $\xi$ is an increasing function with respect to $|x-x_0|$, i.e., the
      limiting term grows for bigger displacements.
\item $\xi(x-x_0, \delta) = 1 \iff |x-x_0| =  \delta$.
      This property is used to give the users a physical intuition of the
      meaning of $\delta$, i.e., the limiting term will start to dominate the
      objective function whenever the mesh nodes move to a distance of $\delta$,
      or bigger, from their original positions.
      Note that, due to their normalization, the metric integrals of
      \eqref{eq_F_full} are always expected to be in $[0,1]$.
\item $\xi$ has two well-defined derivatives with respect to $\bx$.
      This enables us to use a standard Newton-based nonlinear solver.
\item $\xi$ is unitless, i.e., it is invariant to
      rescaling of the domain or change of the physics units system.
\end{enumerate}
The rate at which the limiting term starts to influence the objective function
$F(x)$ depends on the particular form of $\xi$.
Our tests usually utilize one of the following two formulas:
\begin{equation}
\label{eq_xi}
  \xi_1 = \frac{\lvert x - x_0\rvert^2}{\delta^2}, \quad
  \xi_2 = e^{ 10 \left(
              \frac{\lvert x - x_0 \rvert^2}{\delta^2} - 1
              \right) }.
\end{equation}
As shown in Figure \ref{fig_xi}, the exponential option $\xi_2$ becomes active,
and starts to grow much faster, only after mesh displacement approaches the
maximum allowed value of $\delta$.

The constant $c$ in equation \eqref{eq_lim} provides proper normalization.
In particular, note that the integrals in \eqref{eq_lim} are over target
elements, which may or may not have the notion of size.
For example, it is not necessary to set the local volume of a target
element which will be used to perform shape-only optimization.
To preserve invariance under mesh refinement, we set $c$ to
\begin{equation}
\label{eq_lim_scale}
  c =
  \begin{cases}
     \frac{1}{N_E V_{\Omega, avg}}  = \frac{1}{V_{\Omega}}
       & \text{if the targets contain volumetric information}, \\
     \frac{1}{N_E}
       & \text{otherwise},
  \end{cases}
\end{equation}
where $N_E$ is the total number of elements in the mesh, and
$V_{\Omega}$ is the total volume of the domain.

\begin{figure}[h!]
\centerline{
  \includegraphics[width=0.24\textwidth]{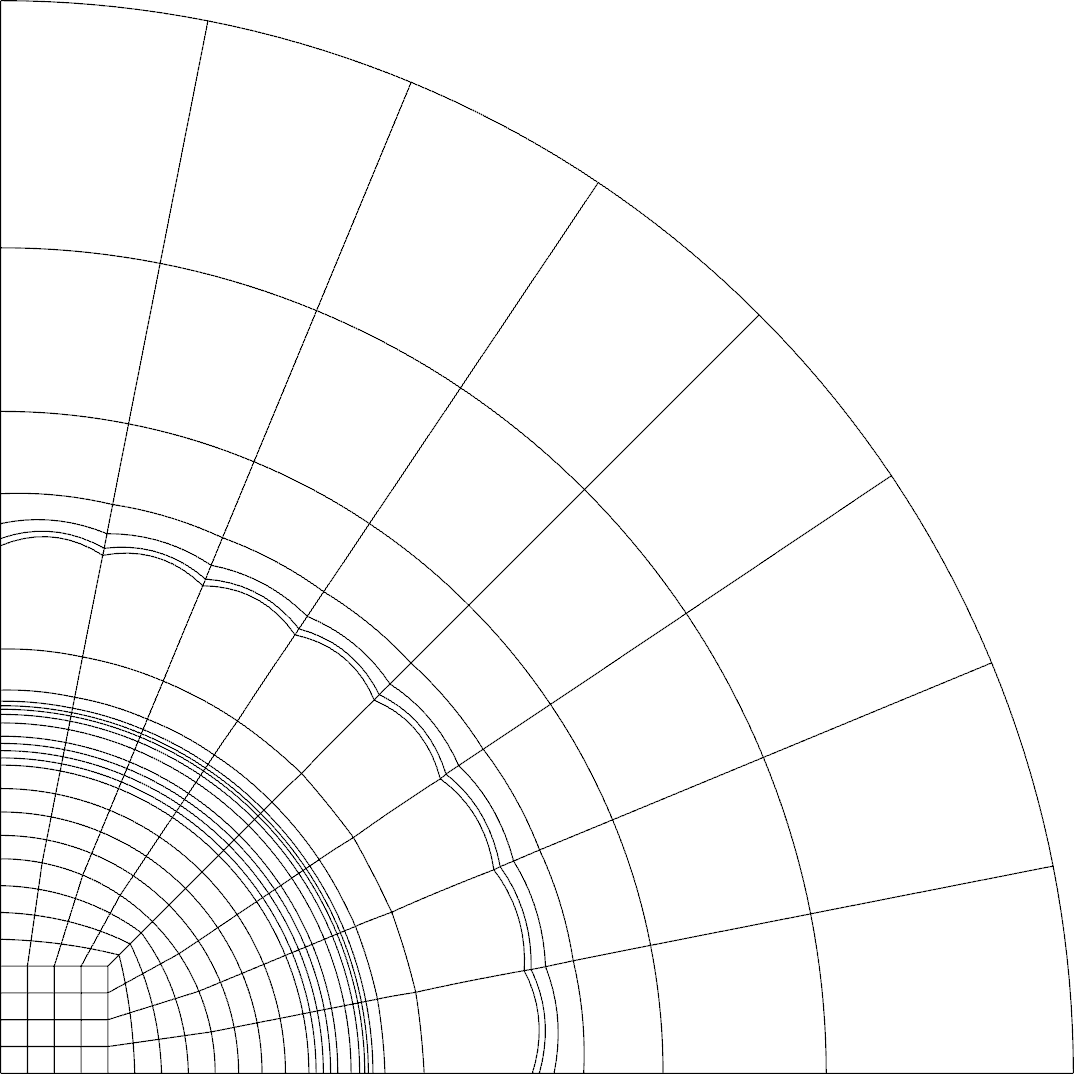} \hfil
  \includegraphics[width=0.24\textwidth]{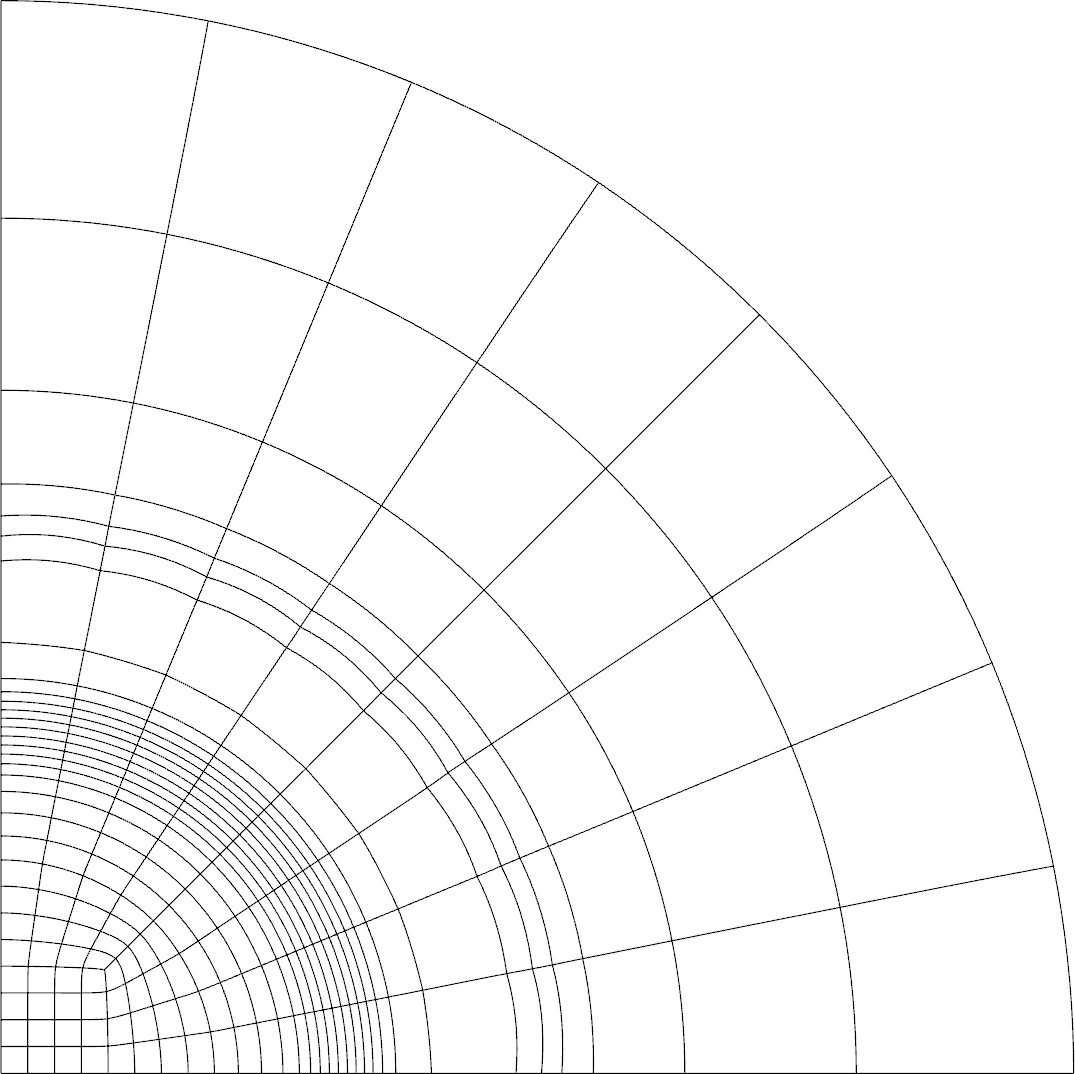} \hfil
  \includegraphics[width=0.24\textwidth]{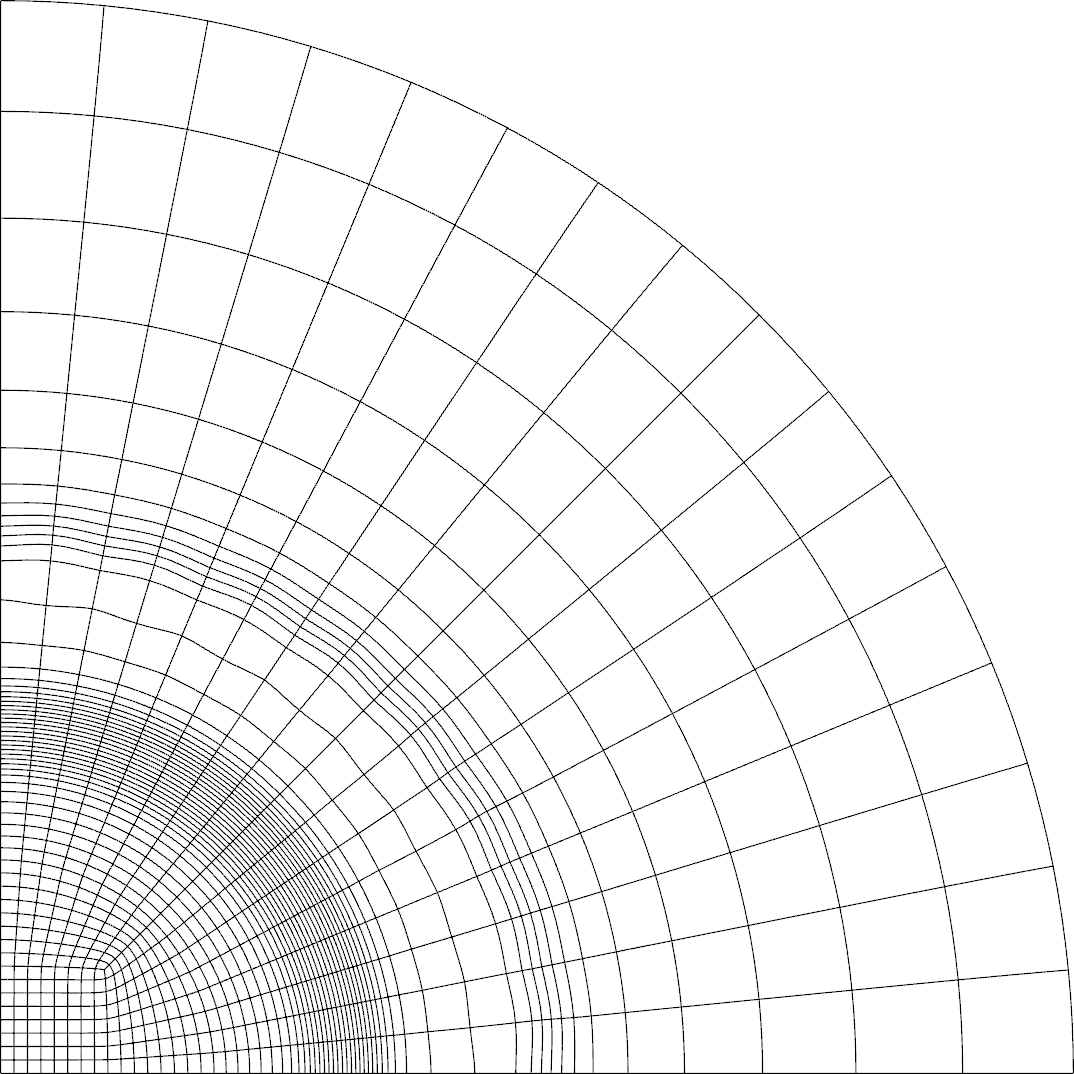} \hfil
  \includegraphics[width=0.24\textwidth]{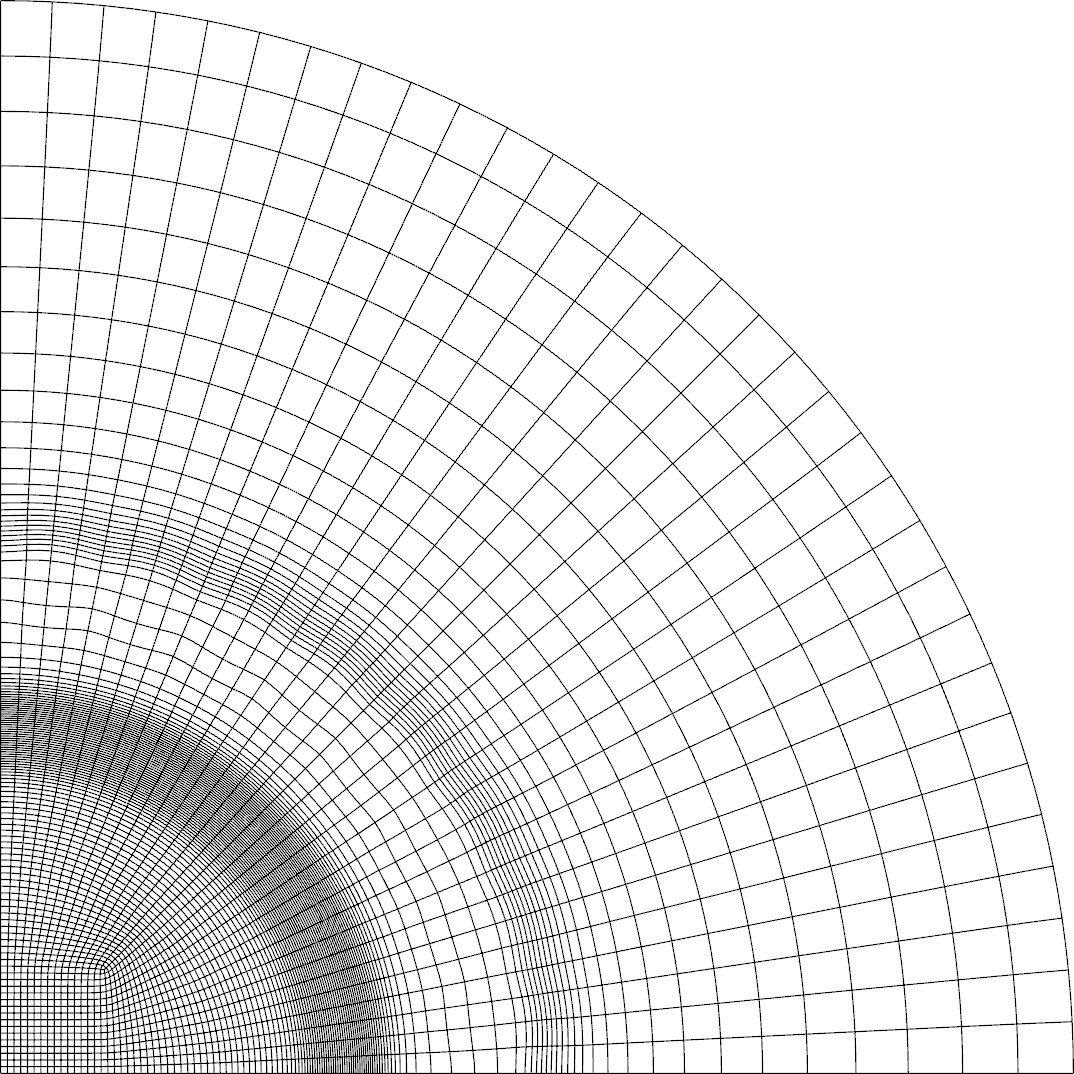}}
\caption{Perturbed initial mesh (left) and optimized meshes
         using $\mu_9$ and $\xi_1$ with $\delta(x_0) = 0.1$
         for 0,1,2 refinements of the original mesh}
\label{fig_icf_normal}
\end{figure}

A simple 2D example that demonstrates the desired behavior of the
normalization factors under mesh refinement is presented in
Figure \ref{fig_icf_normal} and Table \ref{tab_icf_normal}.
In this example we optimize a third order mesh towards ideal shape equally-sized
targets with the shape+size metric $\mu_9=\tau\mid T-T^{-t}\mid^2$.
Node displacements are restricted by using $\xi_1$ and $\delta(x_0) = 0.1$.
We observe similar behavior for the different mesh resolutions.

\begin{table}[h!]
\begin{center}
  \begin{tabular}{c | c c c c}
  \hline
  Refinements & Final $F$ & Metric part & Limiting part & Max displacement  \\
  \hline
  0  & 0.8423 & 0.7819 & 0.0604 & 0.0083  \\
  1  & 0.8422 & 0.7819 & 0.0602 & 0.0083  \\
  2  & 0.8422 & 0.7819 & 0.0602 & 0.0083  \\
  \hline
  \end{tabular}
\end{center}
\caption{Results for final functional value $F$ and its components for several
         refinements of the $Q_3$ ICF mesh,
         using $\mu_9$ and $\xi_1$ with $\delta(x_0) = 0.1$.}
\label{tab_icf_normal}
\end{table}

\begin{figure}[h!]
\centerline{
  \includegraphics[width=0.30\textwidth]{figures/icf_init} \hfil
  \includegraphics[width=0.30\textwidth]{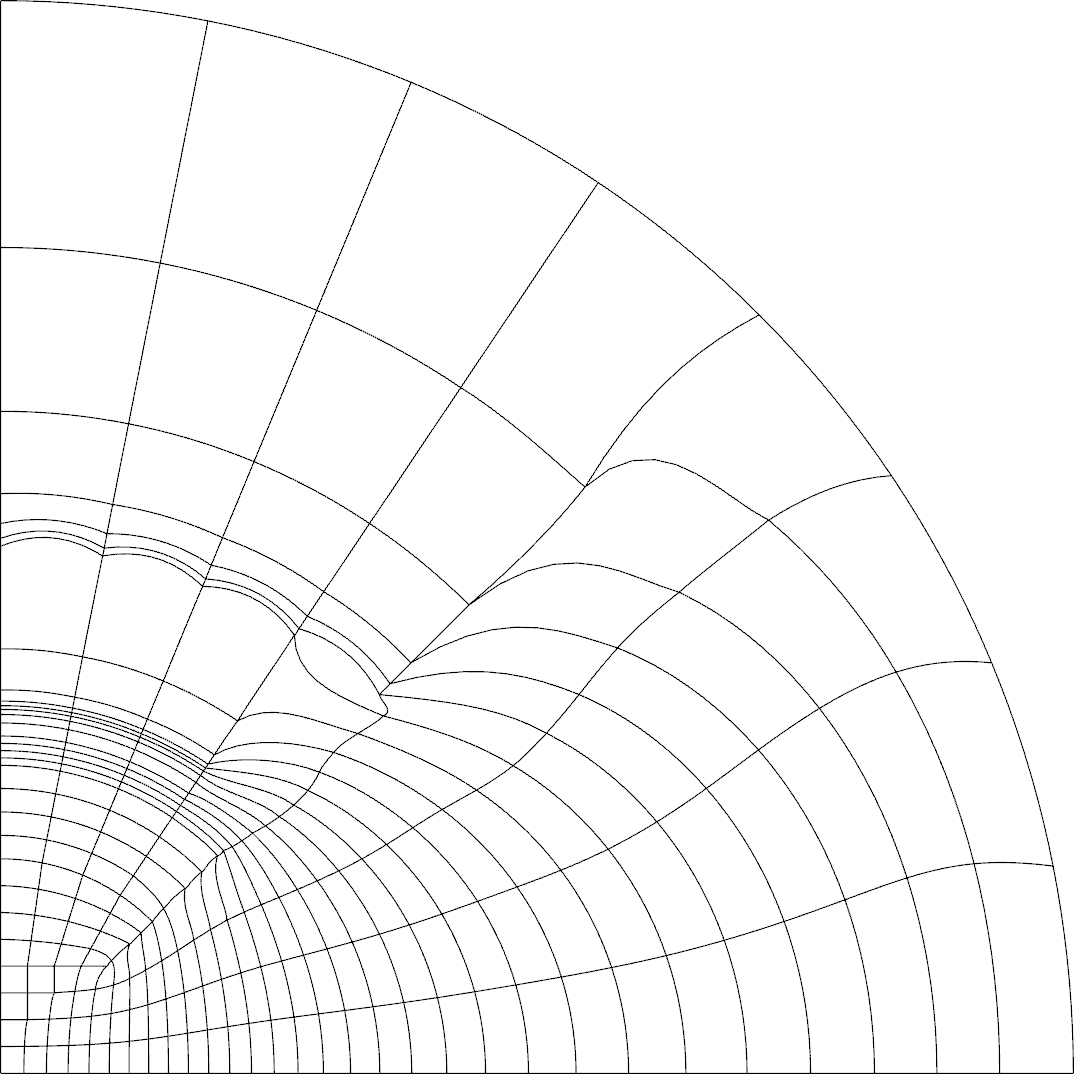} \hfil
  \includegraphics[width=0.30\textwidth]{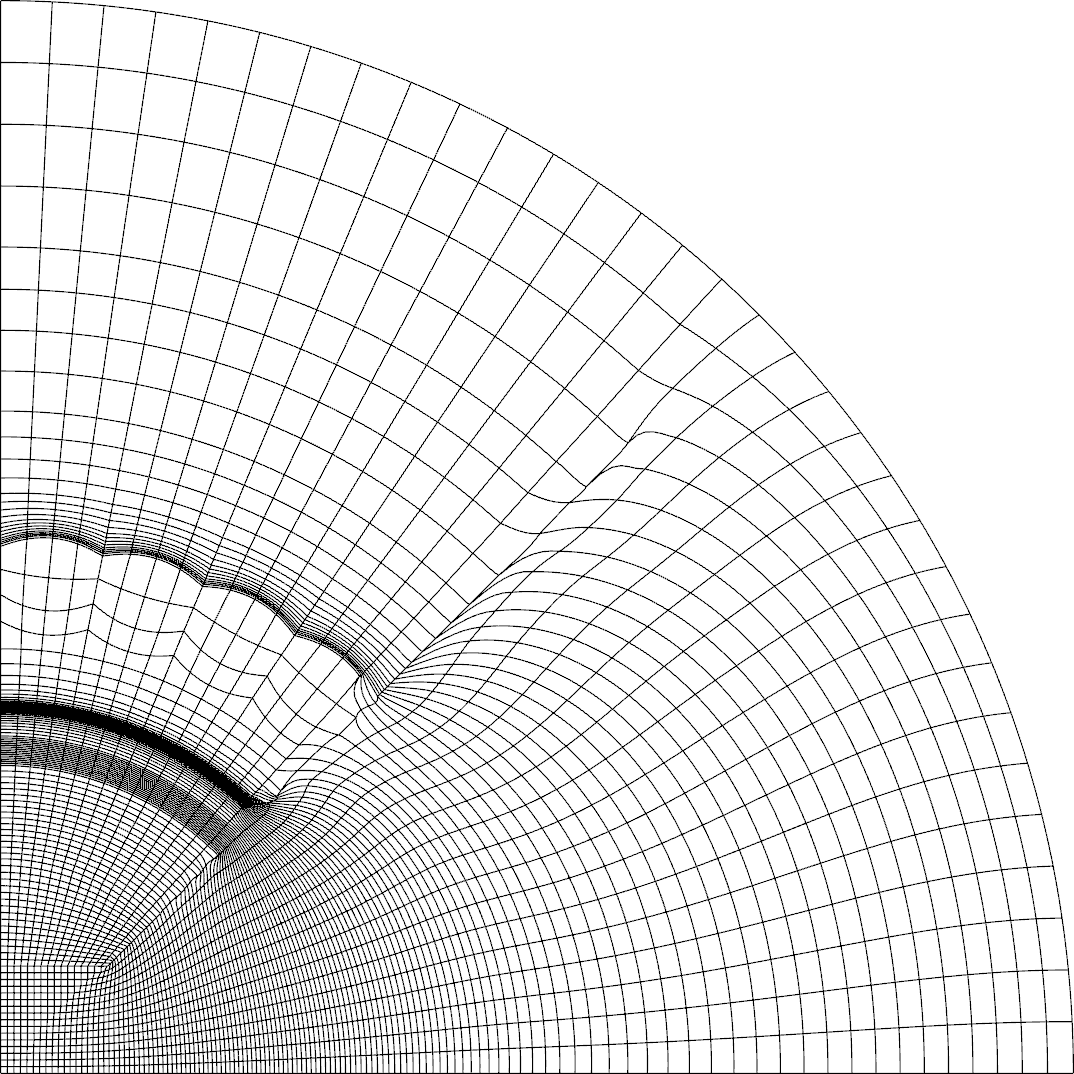}}
\caption{Perturbed initial mesh (left) and locally optimized meshes
         (0 and 2 refinements of the original mesh) using
         space-varying coefficient $\delta(x_0)$. }
\label{fig_icf_loc}
\end{figure}

\paragraph{Local mesh optimization}
A standard approach to perform local mesh optimization is to eliminate unknowns
on the linear algebra level, i.e., to apply the Newton solver on a submesh.
When the region of interest changes dynamically in ALE simulations, however,
refactoring the linear algebra problem at every remesh step can be cumbersome.
Another application of the space-dependent coefficient $\delta(x_0)$
is that it can be used to perform local mesh optimization.
This is demonstrated by another simple 2D test in Figure \ref{fig_icf_loc}.
Again we use a third-order mesh, $\mu_9$, and ideal equally-sized targets,
but in this case we vary the $\delta$ function in space through
\[
  \delta(x_0) =
  \begin{cases}
     10^{-4} & \text{if } x > y \text{~(forcing strong limiting)}, \\
     1.0  & \text{otherwise (effectively no limiting)}.
  \end{cases}
\]
We observe that the chosen metric is active only in the non-restricted region.

\paragraph{Combination with ALE Hydrodynamics}
As already mentioned, the above methods to restrict the node displacements,
or perform local mesh optimization, are useful when one wants to preserve a
certain feature that is available in the Lagrangian mesh.
Our experience shows that it is a good practice to set $\delta(x_0)$ to be some
fraction of the displacement that occurs between the last two ALE steps,
which we denote by $\delta_L(x_0)$.
This ties the amount of mesh relaxation to the amount of Lagrangian motion,
and in particular, guarantees that nodes that do not move during the
Lagrange phase ($\delta_L = 0$), do not move in the remesh phase either.
Various modifications can be performed to relax this restriction in certain
problems, e.g., diffusing $\delta_L(x_0)$ to obtain transition regions.

\section{Mesh Adaptivity and ALE Triggers}
\label{sec_adapt}

Mesh adaptivity \cite{huang2010adaptive} is beneficial when a certain dynamic
feature is not resolved by the Lagrangian mesh.
In this section we demonstrate how the TMOP-based framework adapts the
different geometric properties \eqref{eq_decomp} of the high-order mesh to
control the numerical errors in ALE methods.
We start by describing the coupling between target
construction and discrete simulation data, followed by a discussion of the
appropriate triggering mechanisms in the remesh phase of ALE simulations.


\subsection{Adaptivity to Discrete Simulation Fields}
\label{sec_adapt_discrete}

The objective of the r-adaptivity process is to optimize the mesh using
information from a discrete function, e.g., a finite element solution
function that is evolved during the Lagrangian phase.
In TMOP, r-adaptivity is achieved by incorporating the
discrete data into the target Jacobian matrices.
This process can be split into five major steps:
\begin{enumerate}
  \item Choice of adaptation goal, i.e., which geometric properties of the mesh
        must be controlled and what is their desired behavior.
        This can be motivated, for example, by known a-priori error estimates.
  \item Derivation of geometric parameters from relevant simulation data.
  \item Definition of pointwise target Jacobian matrices $W$.
  \item Choice of a quality metric $\mu(T)$ that measures differences in the
        geometric quantities of interest.
  \item Mesh optimization by solving the final nonlinear optimization problem.
\end{enumerate}
The first two steps are simulation and problem-dependent, but once the
geometric parameters are known, one can construct the targets and optimize the
mesh automatically.
Alternative approaches for linear meshes in the context of ALE hydrodynamics
are described in \cite{Greene2017, Vachal2011}.

\begin{figure}[tbh]
\begin{center}
$\begin{array}{ccc}
\includegraphics[height=0.3\textwidth]{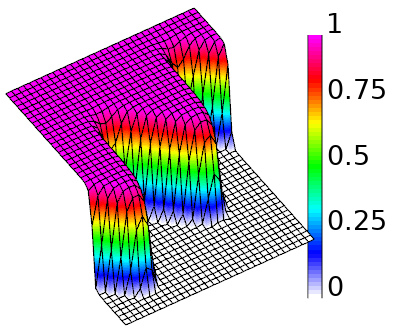} &
\includegraphics[height=0.3\textwidth]{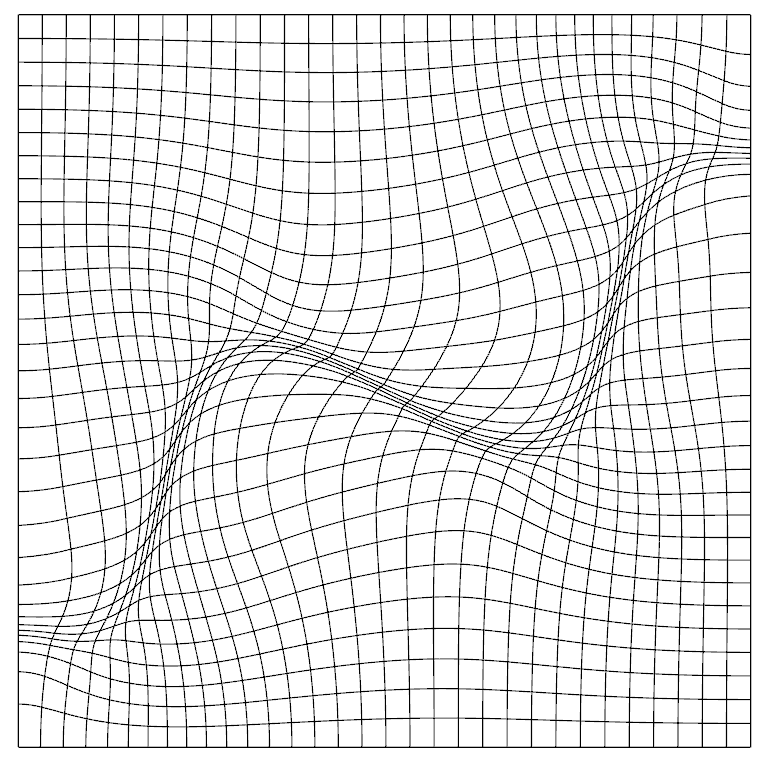} \\
\end{array}$
\end{center}
\caption{Material indicator function on the initial uniform mesh (left panel)
         and adaptively optimized mesh (right panel).}
\label{fig_2matind}
\end{figure}

Next we provide a simple 2D example to illustrate the above procedure.
Suppose we want to adapt the mesh size and aspect ratio to a
material interface, and suppose that the material position is prescribed
by a volume fraction function $\eta$, as it is the case in BLAST,
see Figure \ref{fig_2matind}.
The needed geometric parameters can be computed from the gradient of $\eta$.
The aspect ratio $r$ is computed by the ratio of the gradient components
(i.e., $r \propto \nabla_x\eta/\nabla_y\eta$)
and the local size $s$ depends on the magnitude of the gradient
(i.e., $s \propto ||\nabla\eta||$).
To complete the definition of \emph{shape}, we also choose the target skewness
to be the same as that for an ideal element
($\phi=\pi/2$ for a quadrilateral element).
Using this approach, we construct the target matrix $W$
corresponding to \eqref{eq_decomp} as
\begin{eqnarray}
\label{eq_sinewavetarget}
W =
\begin{bmatrix}
\sqrt{s} & 0 \\
0 & \sqrt{s}
\end{bmatrix}
\begin{bmatrix}
1 & 0 \\
0 & 1
\end{bmatrix}
\begin{bmatrix}
1 & \cos\,\phi \\
0 & \sin\,\phi
\end{bmatrix}
\begin{bmatrix}
\frac{1}{\sqrt{r}} & 0 \\
0 & \sqrt{r}
\end{bmatrix}.
\end{eqnarray}
Next we choose the shape+size metric $\mu_9$, which measures differences of
size, aspect ratio, and skewness between target and physical configurations.
This metric is invariant to orientation (note that $W$ sets targets that have no
rotation, but this is irrelevant due to the choice of $\mu_9$).
The final result for a third order mesh is presented in
Figure \ref{fig_2matind} which shows that the mesh resolves the interface.
Further details about target construction algorithms
can be found in \cite{knupp2019target}.
Note that the above procedure can be applied to various kinds of
simulation data coming from the Lagrangian phase, e.g.,
shock positions, high temperatures, pressure gradients, etc.
Alternative adaptivity approach, for the case when $\eta$ is prescribed
analytically, is presented in \cite{Roca2018}.
We also stress the work presented in \cite{marcon2019rp}, where the
authors utilize an approach similar to TMOP to adapt curved meshes in the
context of steady-state Navier Stokes equations.

\paragraph{Discrete data on intermediate meshes}
Since the discrete solution fields are available only on the
Lagrangian mesh $\mathcal{M}_0$, adaptive mesh optimization requires a method
to evaluate these fields on the intermediate meshes that are obtained during
the Newton iteration.
All simulations in this paper use the advection-based transfer as
described in Section 4.2 of \cite{IMR2018}.
Another option is high-order interpolation between meshes,
see Section 2.3 of \cite{mittal2019nonconforming} and Section 4.1 of
\cite{IMR2018}.

\paragraph{Gradients of adaptive target matrices}
As we use gradient-based nonlinear solvers to update mesh positions, we have to
compute derivatives of $\mu(T=AW^{-1})$ with respect to $\mathbf{x}$.
Since $W$ depends on a discrete function of the mesh positions, the derivatives
$\partial W / \partial \mathbf{x}$ must be taken into account.
The formulas for these derivatives and further discussion are
presented in Section 3 of \cite{IMR2018}.

\paragraph{Combination of adaptivity and limiting of displacements}
Since mesh adaptivity may require significant amounts of mesh displacements to
achieve the desired geometric characteristics, including the limiting term \eqref{eq_lim} in the nonlinear functional can be counterproductive.
Nevertheless, the limiting effects (e.g. local mesh optimization) can still be
achieved by designing a proper space-dependent limiting function $\delta(x_0)$,
which must be aware of the adaptivity goals.


\subsection{Adaptive ALE Triggers}
\label{sec_trigger}

We have briefly discussed ALE triggers in Section 5 of \cite{IMR2018}.
In this section we make several important improvements to the work presented
in \cite{IMR2018}, along with a more detailed discussion.

Since BLAST can perform any number of Lagrangian steps between two ALE steps,
an important question is when, or how often, to optimize the mesh.
Ideally we want a trigger mechanism that does remesh often enough to avoid
small computational time steps due to mesh distortion.
At the same time we do not want to remesh too often,
as ALE calculations can lead to increased numerical dissipation.
In BLAST, this dissipation comes from transitions between
primal and conservative variables and monotonicity treatment
in the remap phase \cite{Dobrev2018}, which is a purely numerical process
that is independent of the physics of the Euler equations.
Furthermore, in the case of mesh adaptivity, higher number of remesh steps
slow down the overall simulations as the mesh adaptivity process involves
the extra computations explained in Section \ref{sec_adapt_discrete}.

The above discussion suggests that the ALE trigger mechanism must be aware of
the mesh quality during the Lagrangian phase, and induce an ALE step
only when the mesh quality is below a certain threshold.
It is important to note that \emph{mesh quality} must be viewed in terms
of the mesh quality metric used in the mesh optimization phase;
the mesh optimization phase is expected to produce a mesh that is above the
trigger threshold.

The TMOP concepts can be utilized to define such mesh optimization triggers,
both in the purely geometric and the simulation-driven case.
Similar to constructing target-matrices, the user specifies
one or more \textit{admissible} Jacobian matrix, $S$, which defines
the transformation from the reference element to the worst element
that can be used during the Lagrangian phase.
Note that $S$ is initialized in the beginning and stays
constant throughout the simulation.
Then $U = S W^{-1}$ represents the Jacobian of the transformation between
target and admissible coordinates, which is used to calculate $\mu(U)$, the
highest admissible mesh quality metric value for the metric of interest.
With that, remesh is triggered whenever $\mu(T) \geq \mu(U)$ at
any mesh quadrature point.
In the case of mesh adaptivity, the target matrices $W$ adapt in time
as they depend on some dynamic solution fields.
Thus $U$ and $\mu(U)$ change too, leading to a trigger that is adapted to the
current Lagrangian solution.

For example, suppose we adapt small mesh size to high temperature regions.
Then at a fixed point in space $x$, the admissible metric value
$\mu(U(x))$ represents the \emph{admissible} size deviation,
given the \emph{current temperature} at $x$.
The value $\mu(T(x))$ represents the \emph{actual} size deviation of the
Lagrangian mesh, at the current temperature.

The specific formulation of $S$ is, of course, problem dependent, e.g., it can
be used to trigger remesh steps whenever particular mesh configurations occur.
Just as the target-matrices $W$, the admissible matrices can contain information
about the local shape, size or alignment, at any quadrature point.
As a simple 2D example, we can define $S$ as
\begin{equation}
\label{eq_U}
S = \begin{pmatrix}
      1 & 0 \\
      0 & S_{22}
    \end{pmatrix} \,,
\text{where } S_{22} > 0 \text{ is a user-specified parameter} \,.
\end{equation}
When $\mu$ is a shape metric, remesh is triggered based on local aspect ratio.
When $\mu$ is a shape+size metric, remesh is triggered based on local size and
aspect ratio.

ALE triggers can be combined with the use of the limiting term \eqref{eq_lim},
as long as the limiting distances $\delta(x_0)$ are large enough to allow the
optimization process to achieve $\mu(T) < \mu(U)$ at all quadrature points.


\section{Numerical tests}
\label{sec_tests}

In this section we report numerical results from the algorithms in the previous
sections as implemented in the BLAST code \cite{blast}.
The baseline TMOP methodology is part of the MFEM finite element library
\cite{mfem, MFEMPaper2019}, which is freely available at \url{mfem.org}.
The application-specific target constructions, however, are only included in
the BLAST implementation.
All physical units in the following tests are based on the
$(cm,g,\mu s)$ unit system.

In what follows, tests that adapt mesh size always utilize a size indicator
function $g(x)$ with values in $[0,1]$, where values close to $1$ represent
regions with smaller mesh size.
Target matrices are constructed as
\begin{equation}
\label{eq_W_size}
  W(x) = [g(x) s + (1 - g(x) ) \alpha s ]^{1/d} I,
\end{equation}
where $\alpha$ is a user-specified size ratio between big and small local size.
Unless otherwise specified, we use $\alpha=10$ for all numerical tests.
Such targets specify custom size, ideal shape, and zero rotation;
in all tests we use metrics that are invariant of rotation.
The small local size $s$ is approximated by taking into account the
total volume $V$ of the domain, the volume $V_g$ occupied by $g(x)$,
the number $N_E$ of available elements, and the specified $\alpha$:
\begin{equation}
\label{eq_size}
  \frac{V_g}{s} + \frac{V - V_g}{\alpha s} = N_E \,,
  \quad \text{where } \quad
  V_g = \int_{\mathcal{M}_0} g_0(x_0) \,, ~~
  V = \int_{\mathcal{M}_0} 1 \,.
\end{equation}
The calculation of $s$ is performed once on the Lagrangian mesh $\mathcal{M}_0$.


\subsection{Triple Point}
\label{sec_3point}

Our first test is the 3-material Triple Point problem
described in \cite{Dobrev2012}.
This test is used to compare different remesh strategies and demonstrate their
advantages and drawbacks in terms of accuracy and computational performance.

\begin{figure}[h!]
\centerline{
  \includegraphics[width=0.7\textwidth]{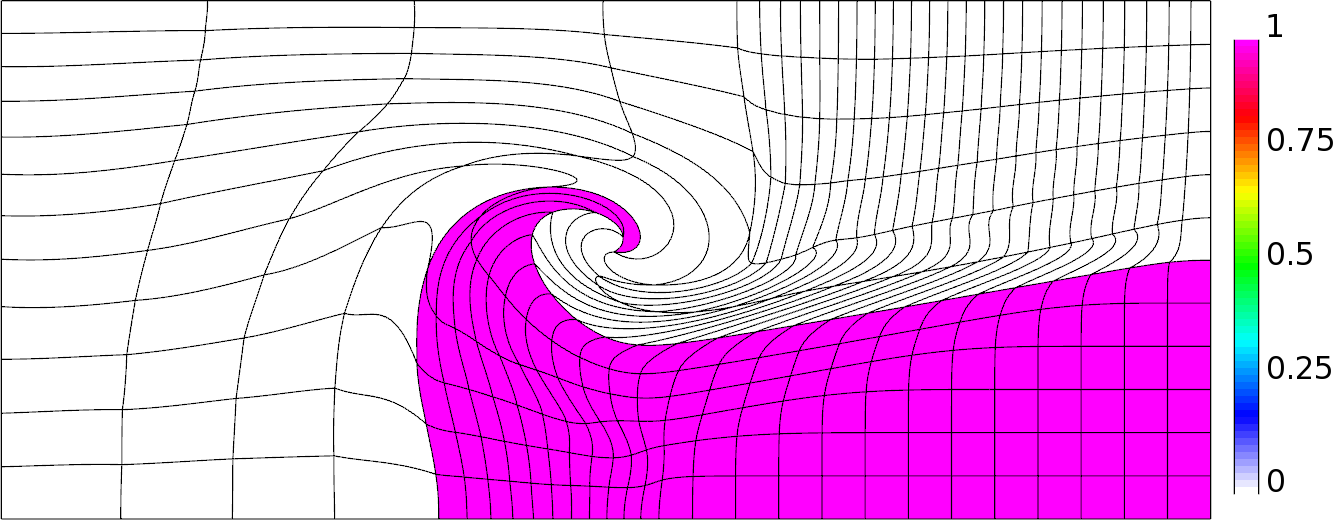}}
\caption{Mesh and material position resulting from a purely
         Lagrangian simulation.}
\label{fig_3p_lagr}
\end{figure}

We start by performing a purely Lagrangian simulation to time 3.5 and record
the final position of one of the materials, see Figure \ref{fig_3p_lagr}.
We take this Lagrangian position as the \emph{true} solution,
as there is no diffusion of the material interface.
For every other test, we define \emph{error} by measuring differences with the
Lagrangian result, for the same material's position.
These differences are measured by sampling both solutions at $10000$
equidistant points over the domain.
This sampling is a nontrivial procedure by itself, as it involves
interpolation in physical space of high-order finite element functions on
curved meshes; for this we use the methods from
Section 2.3 of \cite{mittal2019nonconforming}.
The reported error is the average difference over all points.
All simulations utilize a $Q_3Q_2$ discretization s.t. the mesh is third order.

\begin{figure}[h!]
\centerline{
  \includegraphics[width=0.32\textwidth]{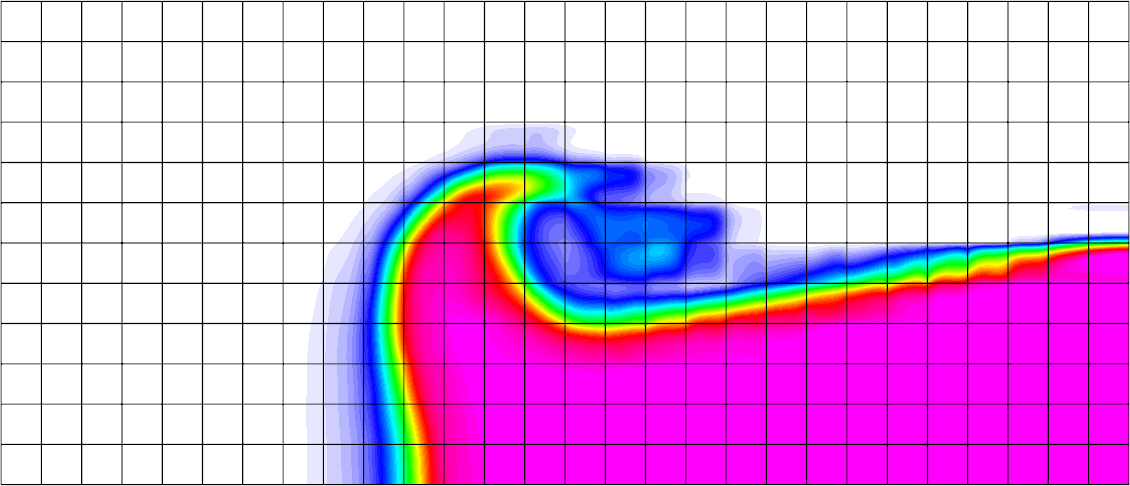} \hfil
  \includegraphics[width=0.32\textwidth]{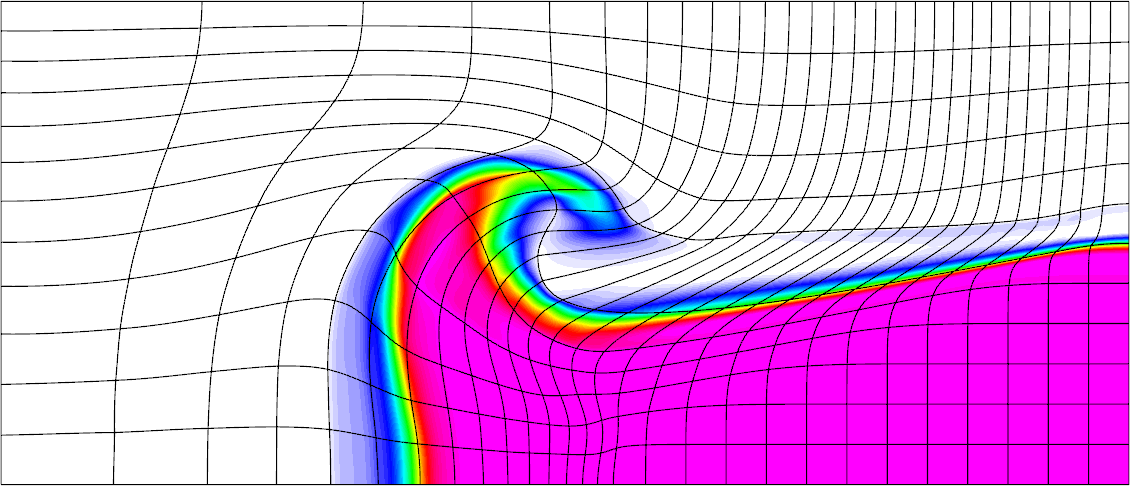} \hfil
  \includegraphics[width=0.32\textwidth]{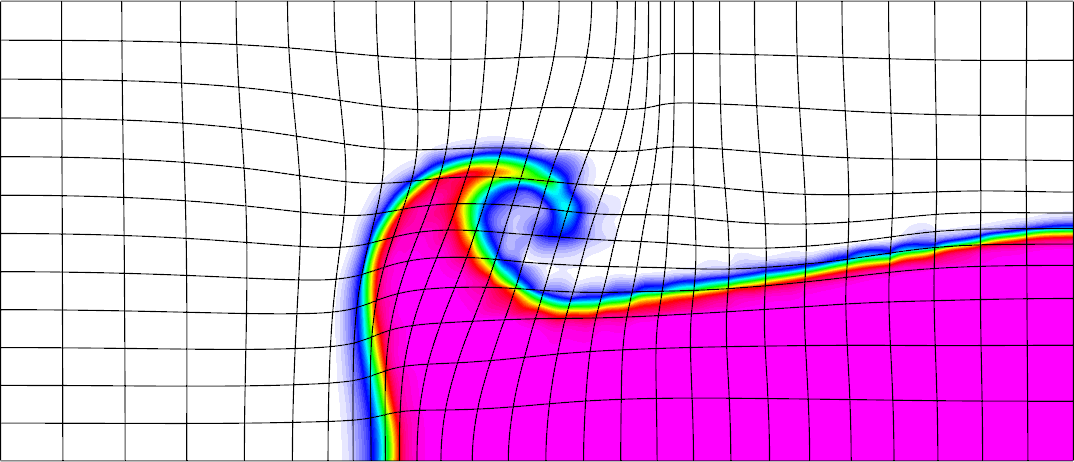}}
\caption{Meshes and material positions resulting from the
         Eulerian (left panel),
         limited by $\delta(x_0)~=~1/2 \delta_L(x_0)$ (center panel), and
         adapted with automatic trigger (right panel) simulations.}
\label{fig_3p_ALE}
\end{figure}

Table \ref{tab_3p} lists all tested approaches.
Most cases use a fixed frequency, so that an ALE step is performed every
50 Lagrangian steps.
The \emph{limited} tests optimize the mesh using the shape metric $\mu_2$
combined with the quadratic $\xi_1$ limiting function and
$\delta(x_0) = a \delta_L(x_0)$, where $a = 1/2$ or $a = 1/3$,
see Section \ref{sec_limiting}.
The \emph{adapted} tests optimize the mesh using the shape+size metric $\mu_7$,
so that material interface regions are assigned smaller mesh sizes.
We report the total number of Lagrangian steps, ALE steps,
remap advection steps (i.e., the total from within all ALE steps),
the relative runtime (all tests use the same computer configuration),
and the error in the material position due to numerical diffusion.
The final material positions for several of the major options is
shown in Figure \ref{fig_3p_ALE}.

\begin{table}[h!]
\begin{center}
  \begin{tabular}{c | c c c c c}
  \hline
  Approach & \# Lag & \# ALE & \# advect & runtime & error \\
  \hline
  Lagrangian                          & 20808 & 0 & 0   & 615 & 0  \\
  Eulerian, period of 50              & 757  & 16 & 90  & 69  & 0.033 \\
  Limited $\frac{1}{2}$, period of 50 & 980  & 20 & 59  & 64  & 0.022 \\
  Limited $\frac{1}{3}$, period of 50 & 1505 & 31 & 55  & 86  & 0.020 \\
  Adapted, auto trigger               & 1666 & 8  & 137 & 174 & 0.017 \\
  Adapted, period of 50               & 1818 & 37 & 154 & 235 & 0.022 \\
  \hline
  \end{tabular}
\end{center}
\caption{Comparison between different ALE strategies for the Triple Point test.}
\label{tab_3p}
\end{table}

We observe that the Eulerian approach, which always returns the mesh to its
original configuration, is fast as it produces large Lagrangian time steps,
but is also the most diffusive.
Performing limiting with $a = 1/2$ is a clear improvement over the
Eulerian case, both in computational time and accuracy, because less advection
remap steps are performed, which compensates the
smaller size of the Lagrangian time steps.
Decreasing the limiting distances by setting $a = 1/3$ preserves the Lagrangian
mesh more closely, resulting in smaller time steps and higher compute time.
The error is slightly better, as less advection steps are performed.

Since the error in this test is always in the interface region, adapting the
mesh size (with adaptive trigger) to this region is the most accurate approach.
The automatic trigger is setup as in \eqref{eq_U} with $S_{22} = 4$, taking
into account the local size and shape of any quadrature point in the mesh.
For example, an ALE would be triggered whenever the local size differs by
a factor of $\geq 4$ of the required mesh size, or the aspect ratio differs by
a factor of $\geq 4$ of the ideal aspect ratio.
Mesh adaptivity, however, leads to slower computations as it leads to bigger
mesh displacements in the remesh phase, which leads to more advection steps
in the remap phase.

As discussed in Section \ref{sec_trigger} performing mesh adaptivity with too
frequent ALE steps can lead to both slower and more diffusive simulations.
This is demonstrated by the last case in Table \ref{tab_3p}, where an ALE
step is performed every 50 steps, independent of the mesh quality during
the Lagrangian phase.


\subsection{2D ALE Simulation of Gas Impact}
\label{sec_gas}

This test represents a high velocity impact of gasses that was originally
proposed in \cite{Barlow14}.
It involves three materials that represent an
{\em impactor}, a {\em wall}, and the {\em background}.
This problem is used to demonstrate the method's behavior in
an impact simulation that cannot be executed in Lagrangian frame
to final time as it produces large mesh deformations.

The domain is $[0, 2] \times [0, 2]$ with $v \cdot n = 0$ boundary conditions.
The material regions are $0.15 \leq x \leq 0.65$ and $0.9 \leq y \leq 1.1$ for
the {\em impactor}, $0.1 \leq x \leq 1.0$ for the {\em wall},
and the rest is {\em background}.
The initial horizontal velocity of the impactor is $0.2 cm / \mu s$.
The problem is run to a final time of $t = 10$ on a $80 \times 80$
second order structured mesh.
The complete thermodynamic setup of this problem and additional
details about our multi-material finite element discretization and overall ALE
method can be found in \cite{Dobrev2016, Dobrev2018}.

The goal of every ALE remesh step is to adapt the size of the
mesh towards the locations of the {\em impactor}, the {\em wall},
and all material interfaces.
The adaptivity field $g(x)$ is constructed from the
simulation's discrete volume fractions $\eta_{wall}$ and $\eta_{imp}$, and
a reconstructed interface indicator function $\eta_{int}$:
\[
g(x) = \max(\eta_{wall}(x), \eta_{imp}(x), \eta_{int}(x)) \,,
\]
where the range of all $\eta$ functions is $[0,1]$.
Target Jacobians $W(x)$ are constructed by \eqref{eq_W_size},
and the shape+size metric $\mu_7$ is used to optimize the mesh.
An automatic ALE trigger is utilized by setting $S_{22} = 4$ in \eqref{eq_U}.

\begin{figure}[!tb]
\centerline{
  \includegraphics[width=0.32\textwidth]{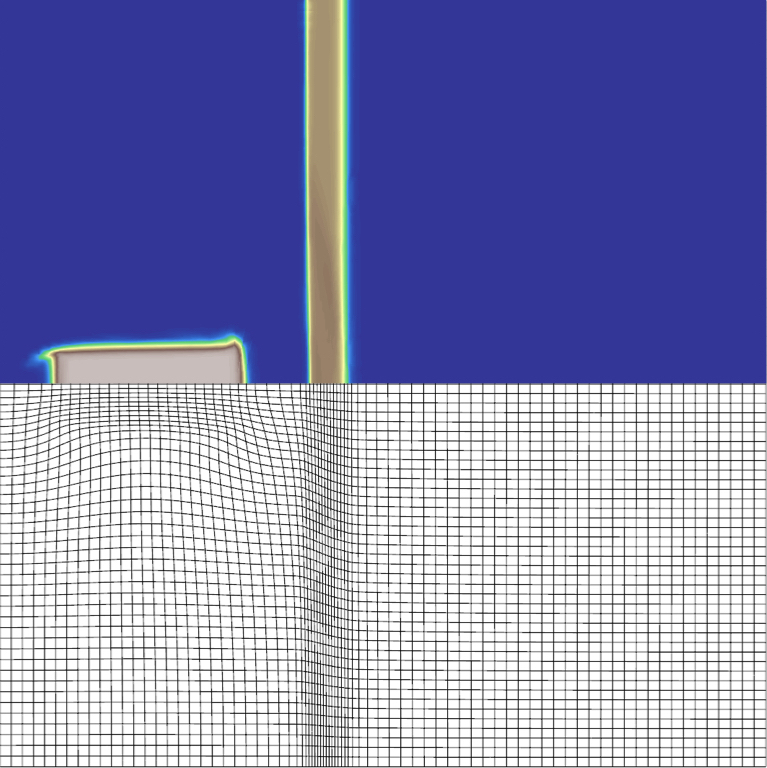} \hspace{0.3mm}
  \includegraphics[width=0.32\textwidth]{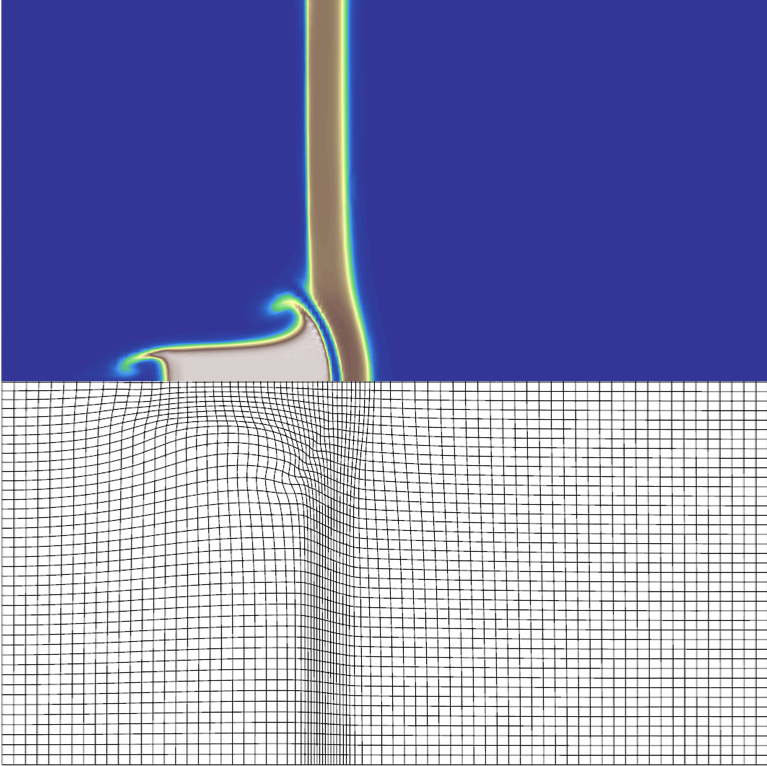}  \hspace{0.3mm}
  \includegraphics[width=0.32\textwidth]{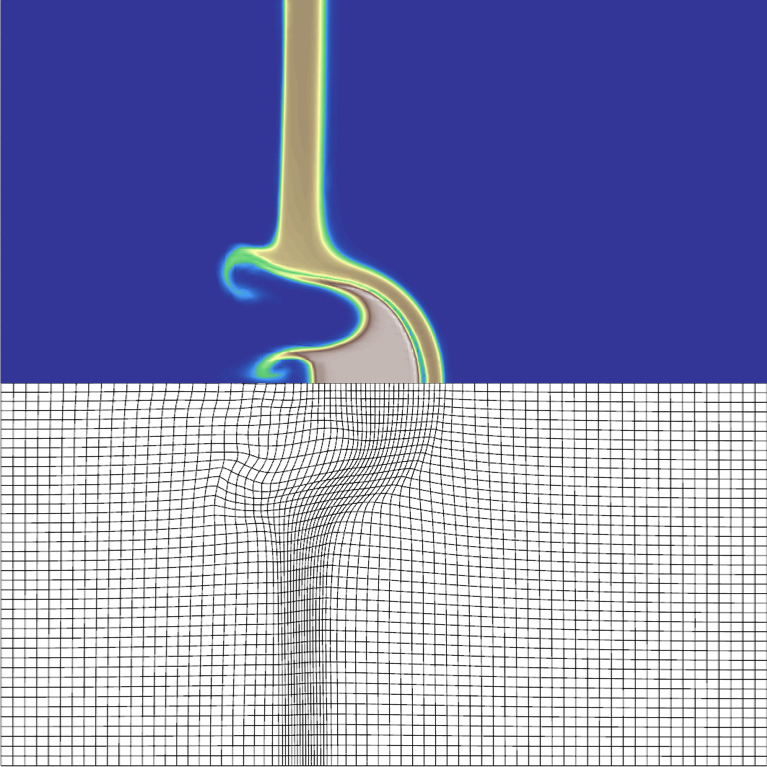}}
\vspace{1.0mm}
\centerline{
  \includegraphics[width=0.32\textwidth]{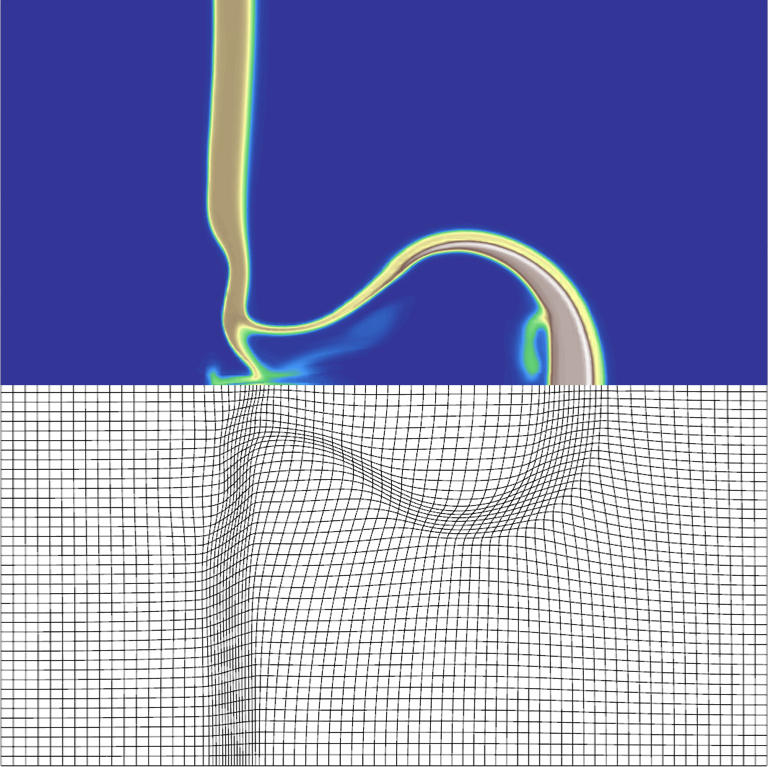}  \hspace{0.3mm}
  \includegraphics[width=0.32\textwidth]{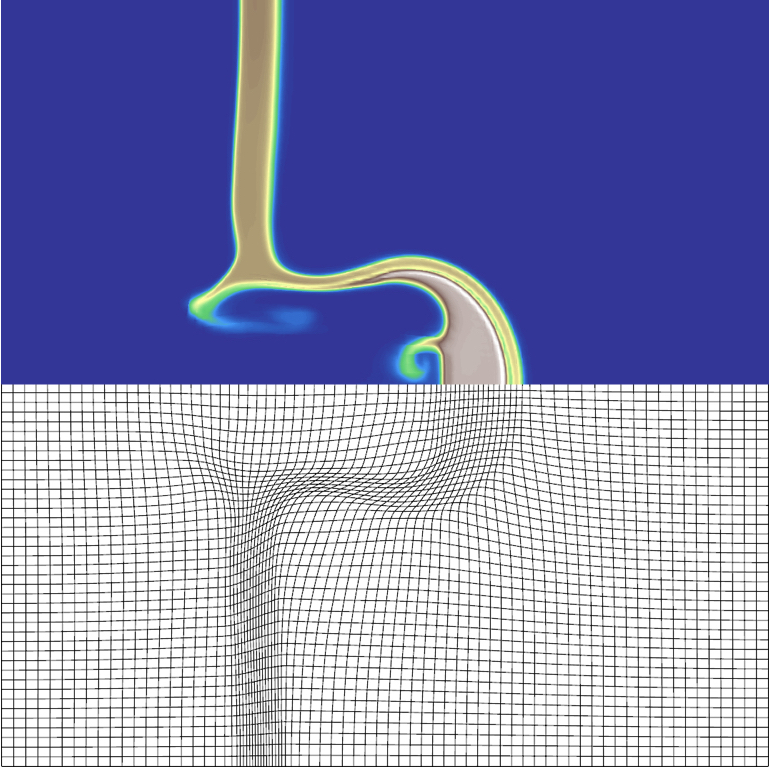}  \hspace{0.3mm}
  \includegraphics[width=0.32\textwidth]{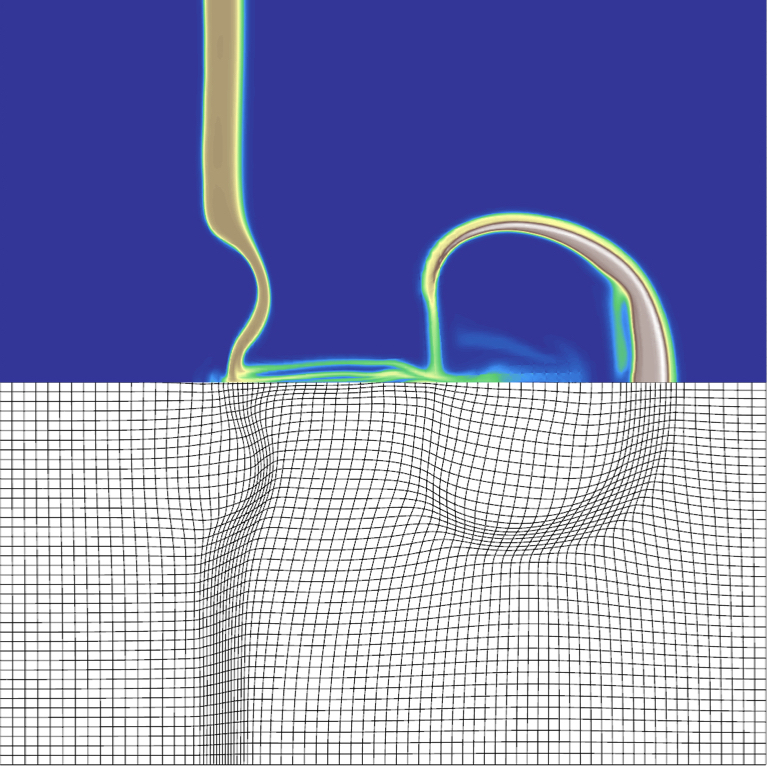}}
\caption{Time evolution of density and mesh position at times
  0.5 (top left), 2 (top middle), 4 (top right), 6 (bottom left),
  8 (bottom middle) and 10 (bottom right) in the 2D gas impact test case.}
\label{fig_gas_impact}
\end{figure}

The time evolution of the material positions and the corresponding mesh
is demonstrated in Figure \ref{fig_gas_impact}.
This calculation required 49400 Lagrange time steps, 38 ALE steps, and 760
advection remap steps in total.
We observe that the algorithm adapts well to the moving materials and the
interface regions, while preserving good overall shape throughout the domain.

\subsection{2Drz and 3D Full ALE Simulation of Steel Ball Impact on Unstructured Mesh}

As a final example, we consider both 2Drz and full 3D calculations of a high
velocity impact of a steel ball against an aluminum plate as
described in \cite{HowellBall02,Dobrev2016}.
The problem consists of a spherical steel projectile
of radius $5.5$ and initial velocity in the z-direction of $0.31 cm / \mu s$ impacting
a cylindrical plate of aluminum with a radius of $24$ and a thickness of $2.5$.
Both the steel and aluminum materials use a Gruneisen equation of state combined with an
elastic perfectly plastic ``strength'' model as described in \cite{Dobrev2014}.
This test demonstrates the ability of the method to adapt 2D and 3D
unstructured meshes in the context of a practical ALE impact simulation.

\begin{figure}[h!]
    \centerline{
      \includegraphics[width=0.32\textwidth]{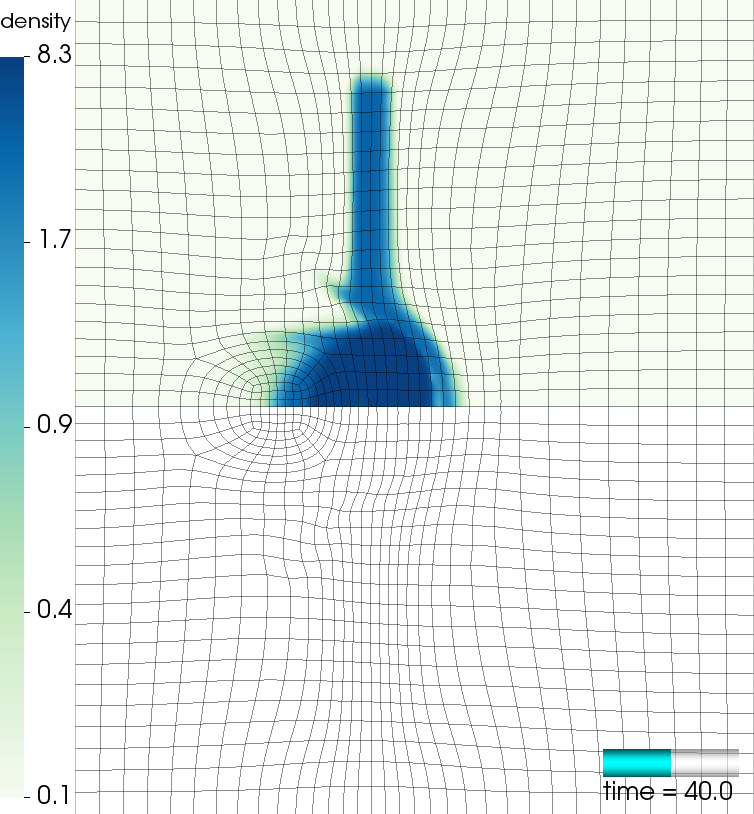} \hfil
      \includegraphics[width=0.32\textwidth]{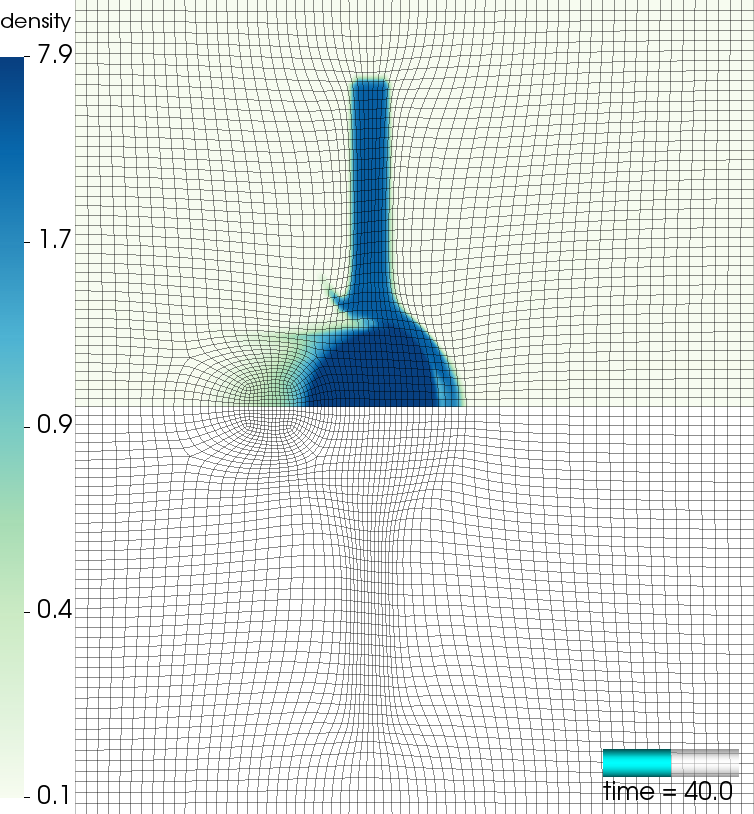} \hfil
      \includegraphics[width=0.32\textwidth]{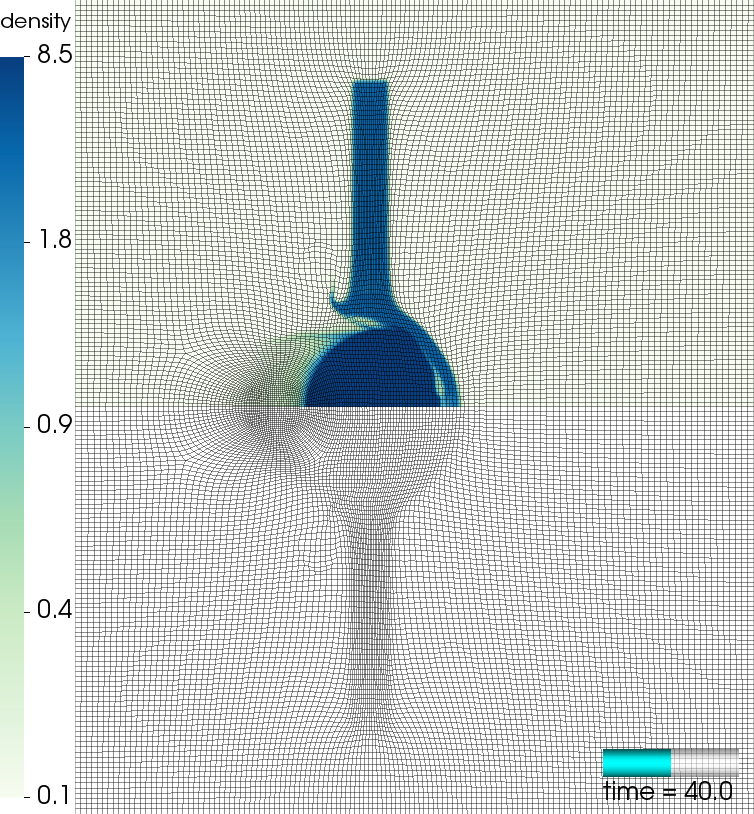}
    }
    \vspace*{4pt}
    \centerline{
      \includegraphics[width=0.32\textwidth]{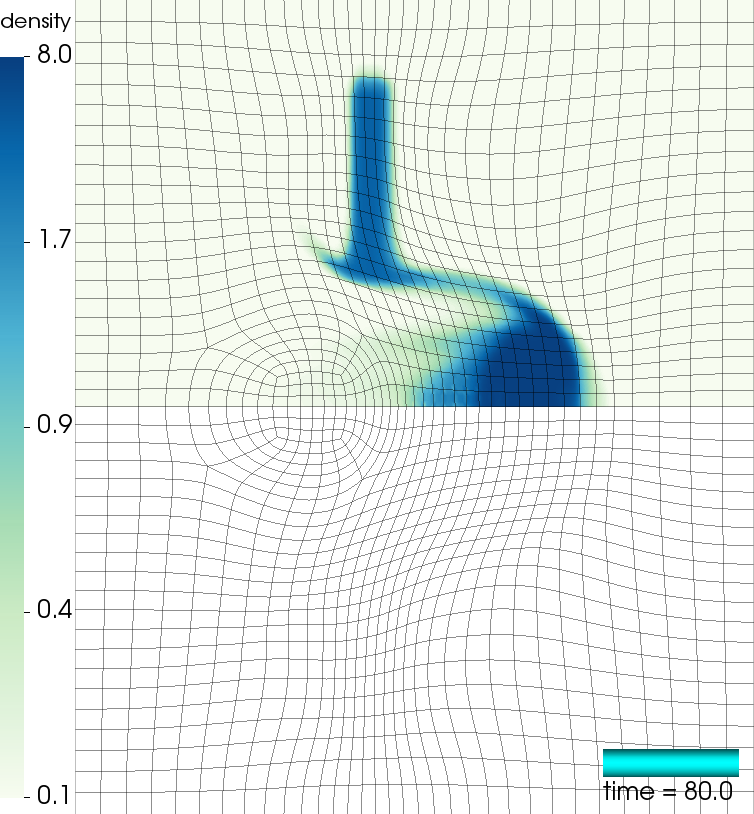} \hfil
      \includegraphics[width=0.32\textwidth]{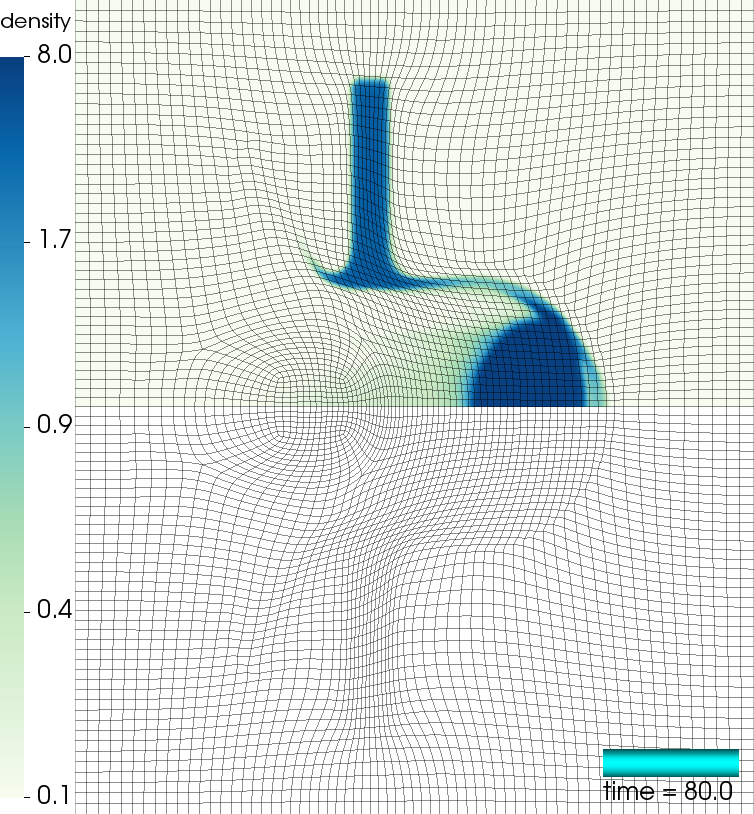} \hfil
      \includegraphics[width=0.32\textwidth]{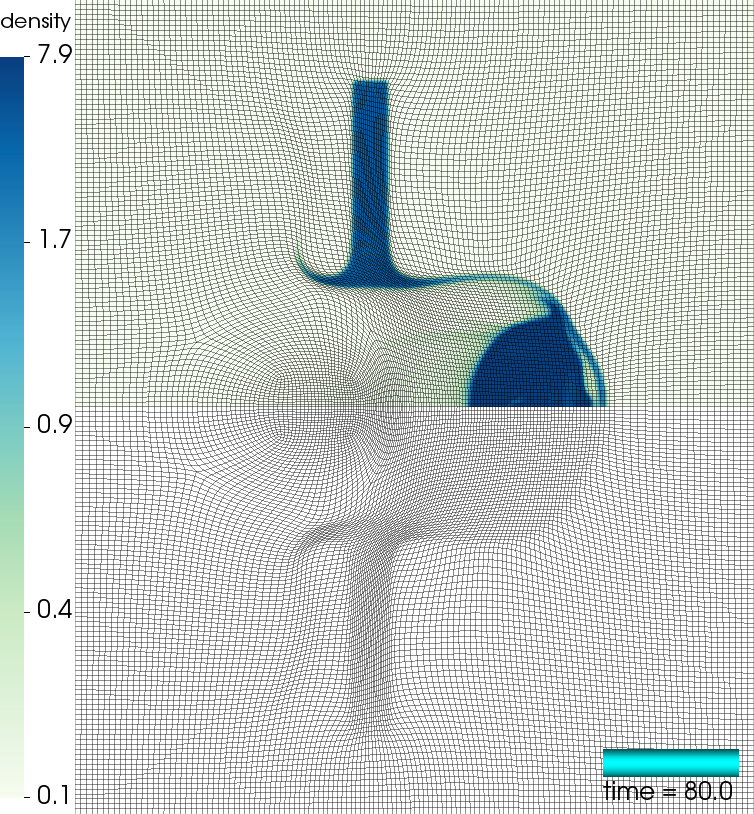}
    }
    \caption{Density and mesh positions at times $t=40$ (\emph{top}) and $t=80$ (\emph{bottom}) using three different mesh
             resolutions for the 2Drz ball impact problem.}
    \label{fig_ball_impact2d}
\end{figure}

In Figure \ref{fig_ball_impact2d} we show results of the 2D (axisymmetric)
simulation on a sequence of refined 2D, high-order, unstructured meshes.
Similarly to the previous example, we adapt to the position of the
{\em ball} impactor, the aluminum {\em plate}, and all material interfaces,
using the shape+size metric $\mu_7$.
The adapted trigger is setup with $S_{22} = 5$ in \eqref{eq_U}.
For the sequence of mesh  refinements, the calculations required a total of
$3201$, $4815$ and $10145$ Lagrange time steps for the coarse, medium and fine
meshes, respectively, to reach a final simulation time of $t=80$.
The total number of remesh steps performed for the each 2Drz calculation is
$19$, $39$ and $88$ for the coarse, medium and fine meshes, respectively.
The total number of remap advection steps is $190$, $382$ and $782$.

\begin{figure}[h!]
    \centerline{
      \includegraphics[width=0.48\textwidth]{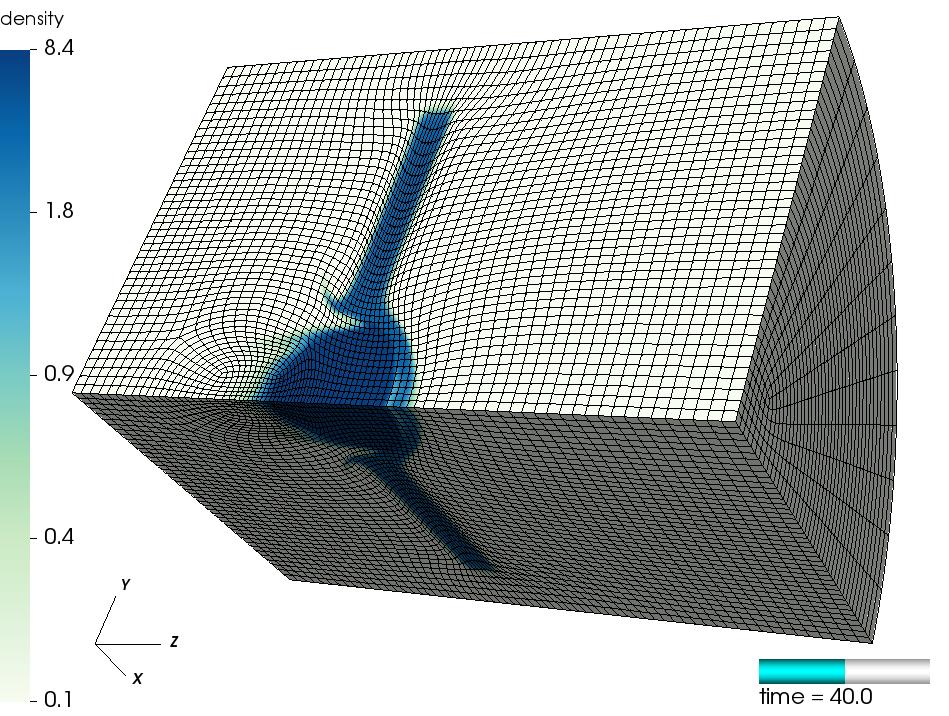} \hfil
      \includegraphics[width=0.48\textwidth]{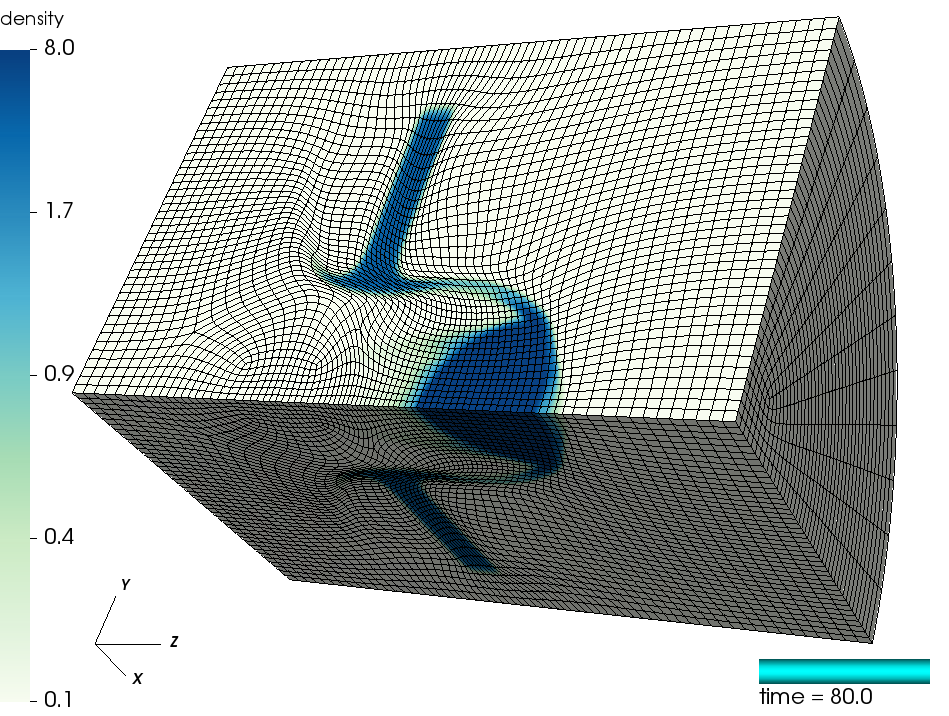}
    }
    \caption{Density and mesh positions at times $t=40$ (\emph{left}) and $t=80$ (\emph{right}) for the 3D ball impact problem.}
    \label{fig_ball_impact3d}
\end{figure}

In Figure \ref{fig_ball_impact3d} we show results of the full 3D simulation on a 
high-order 3D unstructured mesh with resolution identical to the previous 2D
medium case in $rz$ plane.
The adaptivity parameters for this calculation are $\alpha=12$ and
$S_{33} = 8$ for the 3D version of \eqref{eq_U}.
This calculation requires a total of $9671$ Lagrange time steps, $33$ remesh
steps, and $817$ advection remap steps to reach a final simulation time of $t=80$.

\section{Conclusion}
\label{sec_concl}

In this paper we have presented algorithms for simulation-driven optimization
of high-order curved meshes, with application to high-order arbitrary
Lagrangian-Eulerian simulations \cite{blast}.
These allow the user to optimize the shape and size of the mesh,
preserve geometric features, adapt the mesh characteristics to discrete
fields of interest, and control the ALE frequency by adaptive triggers.
These methods provide flexible control over the balance between computational
time and amount of numerical dissipation in ALE simulations.

A future area of research will be to combine the existing TMOP technology
with adaptive mesh refinement.
This would allow to optimize shape, while achieving arbitrary small local size,
which is currently constrained by the fact that the mesh topology is fixed.
We also plan to add capabilities for tangential relaxation around discrete
surfaces, which would be highly beneficial in multi-material ALE simulations.

\noindent \section*{In memoriam}

\noindent This paper is dedicated to the memory of Dr. Douglas Nelson Woods
($^*$January 11\textsuperscript{th} 1985 - $\dagger$September
11\textsuperscript{th} 2019), promising young scientist and post-doctoral
research fellow at Los Alamos National Laboratory.  Our thoughts and wishes go
to his wife Jessica, to his parents Susan and Tom, to his sister Rebecca and to
his brother Chris, whom he left behind.

\noindent \section*{Disclaimer}
This document was prepared as an account of work sponsored by an agency of the
United States government. Neither the United States government nor Lawrence
Livermore National Security, LLC, nor any of their employees makes any warranty,
expressed or implied, or assumes any legal liability or responsibility for the
accuracy, completeness, or usefulness of any information, apparatus, product, or
process disclosed, or represents that its use would not infringe privately owned
rights. Reference herein to any specific commercial product, process, or service
by trade name, trademark, manufacturer, or otherwise does not necessarily
constitute or imply its endorsement, recommendation, or favoring by the United
States government or Lawrence Livermore National Security, LLC. The views and
opinions of authors expressed herein do not necessarily state or reflect those
of the United States government or Lawrence Livermore National Security, LLC,
and shall not be used for advertising or product endorsement purposes.

\bibliographystyle{elsarticle-num}
\bibliography{caf2019}

\end{document}